\documentclass[review]{elsarticle}
\usepackage{lineno,hyperref}
\usepackage{amsmath}
\usepackage{multirow}
\usepackage{mdwlist}
\usepackage{xcolor}
\usepackage{subfigure}
\usepackage{colortbl}
\usepackage{amssymb}
\modulolinenumbers[5]

\journal{Numerical Linear Algebra with Applications}

\biboptions{numbers,sort&compress}

\newtheorem{Theorem}{Theorem}
\newtheorem{Example}{Example}
\newtheorem{Lemma}{Lemma}
\newtheorem{Remark}{Remark}
\newtheorem{Corollary}{Corollary}
\newtheorem{Definition}{Definition}
\newtheorem{Algorithm}{Algorithm}

\newenvironment{Proof}[1][Proof.]{\begin{trivlist}
\item[\hskip \labelsep {\bfseries #1}]}{\ep\end{trivlist}}
\newcommand{\ep}{\hfill\rule{0.15cm}{0.35cm}\vskip 0.3cm}

\begin{document}

\begin{frontmatter}

\title{Time-space Finite Element Adaptive AMG for Multi-term Time Fractional Advection Diffusion Equations}

\author[mymainaddress,mysecondaryaddress]{Xiaoqiang Yue}
\ead{yuexq@xtu.edu.cn}

\author[mymainaddress]{Yehong Xu}

\author[mysecondaryaddress]{Shi Shu\corref{mycorrespondingauthor}}
\cortext[mycorrespondingauthor]{Corresponding author}
\ead{shushi@xtu.edu.cn}

\author[mymainaddress]{Menghuan Liu}

\author[mysecondaryaddress]{Weiping Bu}

\address[mymainaddress]{School of Mathematics and Computational Science, Xiangtan University, Hunan 411105, P.R. China}
\address[mysecondaryaddress]{Hunan Key Laboratory for Computation and Simulation in Science and Engineering,
Xiangtan University, Hunan 411105, P.R. China}

\begin{abstract}
In this study we construct a time-space finite element (FE) scheme and furnish cost-efficient approximations for
one-dimensional multi-term time fractional advection diffusion equations on a bounded domain $\Omega$.
Firstly, a fully discrete scheme is obtained by the linear FE method in both temporal and spatial directions,
and many characterizations on the resulting matrix are established.
Secondly, the condition number estimation is proved,
an adaptive algebraic multigrid (AMG) method is further developed to lessen computational cost and analyzed in the classical framework.
Finally, some numerical experiments are implemented to reach the saturation error order in the $L^2(\Omega)$ norm sense,
and present theoretical confirmations and predictable behaviors of the proposed algorithm.
\end{abstract}

\begin{keyword}
Multi-term time fractional advection diffusion equations,
time-space finite element, condition number estimation, algorithmic complexity, adaptive AMG method
\MSC[2010] 35R11\sep  65F10\sep  65F15\sep  65N55
\end{keyword}

\end{frontmatter}

\linenumbers

\section{Introduction}

In recent years, fractional differential equations (FDEs) have burst onto the scientific computing scene
as a powerful instrument in descriptions of memory and hereditary that has yielded a wide variety of applications
in physics, hydrology, finance and other fields \cite{m-001}. Since the vast majority couldn't be solved accurately,
or their analytical solutions (if luckily derived) always contain specific infinite series resulting in sharp costs of evaluations,
numerical solutions to FDEs becomes very practical and prevalent. Numerous numerical (unconditionally stable and efficient) methods arise,
e.g. finite difference (FD) \cite{l-003,g-001,c-002,w-001,y-002}, finite element (FE) \cite{e-001,z-002,b-002,b-001,m-002,f-001,b-004},
finite volume \cite{l-004} and spectral (element) methods \cite{l-005,z-004,y-005,y-003}.

Fractional advection diffusion equations (FADEs) are known as one of the foremost models in depictions for transport process
in complex systems governed by the abnormal diffusion and non-exponential relaxation patterns \cite{m-003}.
Fundamental and numerical solutions for FADEs with single-term, two-term and multi-term time fractional derivatives
have been investigated in \cite{e-001,s-002,z-005,m-004,b-005,j-001,y-006,b-003}.
Nonetheless, from the survey of references, there are no calculations regarding the time-space FE discretization
for FADEs in literature. In this paper we focus on this topic to
the one-dimensional version of multi-term time fractional advection diffusion equations (MTFADEs)
\begin{align}\label{chp-01-01}
  &\sum_{i=0}^s a_i {}_0^CD_t^{\alpha_i}u(x,t)=K_1\frac{\partial^{2\beta}u(x,t)}{\partial|x|^{2\beta}}
  +K_2\frac{\partial^{2\gamma}u(x,t)}{\partial|x|^{2\gamma}}+f(x,t),~t\in I=(0,T],~x\in \Omega=(a,b) \\
  &\label{chp-01-02}u(a,t)=u(b,t)=0,~t\in I \\
  &\label{chp-01-03}u(x,0)=\psi_0(x),~x\in \overline{\Omega}
\end{align}
with orders $0<\alpha_s<\cdots<\alpha_1<\alpha_0<1$, $\beta\in(0,1/2)$ and
$\gamma\in(1/2,1)$, constants $a_0>0$, $a_i\ge 0$ ($i=1,\cdots,s$), $K_1>0$ and $K_2>0$,
Caputo of order $\alpha$ ($0<\alpha<1$) and Riesz of order $2\varrho$ ($m-1<2\varrho<m$) fractional derivatives defined by
\begin{eqnarray*}
  {}_0^CD^\alpha_tu=\frac{1}{\Gamma(1-\alpha)}\int_0^t(t-s)^{-\alpha}\frac{\partial u}{\partial s}ds,~
  \frac{\partial^{2\varrho}u}{\partial|x|^{2\varrho}}=-\frac{1}{2\cos(\varrho\pi)}({}_xD_L^{2\varrho}u + {}_xD_R^{2\varrho}u)
\end{eqnarray*}
respectively, where
\begin{eqnarray*}
  {}_xD_L^{2\varrho}u = \frac{1}{\Gamma(m-2\varrho)}\frac{\partial^m}{\partial x^m}\int_a^x(x-s)^{m-1-2\varrho}uds,~
  {}_xD_R^{2\varrho}u = \frac{(-1)^m}{\Gamma(m-2\varrho)}\frac{\partial^m}{\partial x^m}\int_x^b(s-x)^{m-1-2\varrho}uds.
\end{eqnarray*}

Another important note is that no matter which fully discrete scheme is utilized, 
there always exists the computational challenge in nonlocality caused by fractional differential operators \cite{g-002}.
Quite many scholars are working to identify algorithms most appropriate to overcome the challenge and utilize computer resources.
Multigrid exploits a hierarchy of grids or multiscale representations, and reduces the error of the approximation
at a number of frequencies (from global smooth to local oscillation) simultaneously \cite{t-001,f-002,f-003,x-001},
which makes multigrid as an extremely superior solver or preconditioner of particular interest.
Pang and Sun presented an efficient V-cycle geometric multigrid (GMG)
with fast Fourier transform (FFT) to solve one-dimensional space-fractional diffusion equations (SFDEs) discretized by an implicit FD scheme \cite{p-001}.
Bu et al. extended and analyzed the V-cycle GMG for one-dimensional MTFADEs via the FD in temporal and FE in spatial directions \cite{b-003}.
Zhou and Wu discussed the FE F-cycle GMG to linear stationary FADEs in Riemann-Liouvlle fractional derivatives \cite{z-003}.
Jiang and Xu constructed optimal GMG approaches for two-dimensional SFDEs to get FE approximations \cite{j-002}.
Chen et al. generalized an algebraic multigrid (AMG) with line smoothers to fractional Laplacian problems
through localizing them as nonuniform elliptic equations \cite{c-003}.
Zhao et al. considered the adaptive FE V-cycle GMG for one-dimensional SFDEs using hierarchical matrices \cite{z-006}.
Chen and Deng exploited GMG's coarsening strategy and grid-transfer operators, equipped with Galerkin coarse-grid operator
to produce a robust multigrid but with much lower convergence rate for nonlocal models with a finite range of interactions \cite{c-001}.
More recently, we developed and analyzed a straightforward adaptive AMG through condition number estimations
for one-dimensional time-space Caputo-Riesz FDEs \cite{y-001}. To the best of our knowledge,
cost-efficient AMG resolutions for MTFADEs by fully time-space FE schemes are still limited.

The goal of this paper is to design a time-space FE scheme and develop a fairly robust and efficient solver
for problem \eqref{chp-01-01}-\eqref{chp-01-03}. The remainder proceeds as follows.
Section 2 contains a review of preliminary knowledge on fractional derivative spaces,
covers the constitution and fundamental properties of the fully discrete FE discretization.
In section 3, condition number estimation on the coefficient matrix is discussed,
followed by the construction and convergence analysis of an adaptive AMG method.
Section 4 reports and analyzes numerical results to showcase the benefits,
and some concluding remarks with follow-up work are given in section 5.

For simplicity, symbols $\lesssim$, $\gtrsim$ and $\simeq$ are used throughout the paper:
$u_1\lesssim v_1$ symbolizes $u_1\leq C_1v_1$, $u_2\gtrsim v_2$ means $u_2\ge c_2v_2$ while $u_3\simeq v_3$ stands for
$c_3v_3\leq u_3\leq C_3v_3$, where $C_1$, $c_2$, $c_3$ and $C_3$ are positive constants independent of step sizes and variables.

\section{Time-space FE scheme for MTFADEs}

In this section we will briefly draw some fractional derivative spaces and several relevant auxiliary results,
which is the basis of our description of the time-space FE scheme in section 2.2,
where also address numerous features of the resulting stiffness and coefficient matrices.

\subsection{Reminder on fractional calculus}

\begin{Definition} (Left and right fractional derivative spaces)
  For any constant $\mu>0$, define norms
\begin{eqnarray*}
  \|u\|_{J_L^\mu(\Omega)} := (\|u\|^2_{L^2(\Omega)}+\|{}_xD_L^\mu u\|^2_{L^2(\Omega)})^{\frac{1}{2}},~
  \|u\|_{J_R^\mu(\Omega)} := (\|u\|^2_{L^2(\Omega)}+\|{}_xD_R^\mu u\|^2_{L^2(\Omega)})^{\frac{1}{2}},
\end{eqnarray*}
  and let $J_{L,0}^\mu(\Omega)$ and $J_{R,0}^\mu(\Omega)$ be closures of $C_0^\infty(\Omega)$
  under $\|\cdot\|_{J_L^\mu(\Omega)}$ and $\|\cdot\|_{J_R^\mu(\Omega)}$, respectively.
\end{Definition}

\begin{Definition}\label{def-01} (Fractional Sobolev space)
  For any constant $\mu>0$, define the norm
\begin{eqnarray*}
  \|u\|_{H^\mu(\Omega)} := (\|u\|^2_{L^2(\Omega)}+\||\xi|^\mu\tilde{u}\|^2_{L^2(\Omega_\xi)})^{\frac{1}{2}},
\end{eqnarray*}
  and let $H_0^\mu(\Omega)$ be the closure of $C_0^\infty(\Omega)$ in $\|\cdot\|_{H^\mu(\Omega)}$ sense,
  where $\tilde{u}$ is the Fourier transform of $u$.
\end{Definition}

\begin{Lemma} (see \cite{z-002}, Proposition 1)\label{clp-02-01}
  If constant $\mu\in(0,1)$ and $u,v\in J_{L,0}^{2\mu}(\Omega)$ (or $J_{R,0}^{2\mu}(\Omega)$), then
\begin{eqnarray*}
  ({}_xD_L^{2\mu}u,v)_{L^2(\Omega)}=({}_xD_L^\mu u,{}_xD_R^\mu v)_{L^2(\Omega)},~
  ({}_xD_R^{2\mu}u,v)_{L^2(\Omega)}=({}_xD_R^\mu u,{}_xD_L^\mu v)_{L^2(\Omega)}.
\end{eqnarray*}
\end{Lemma}

\begin{Lemma} (see \cite{e-001}, Lemma 2.4)\label{clp-02-02}
  For any constant $\mu > 0$ and real valued function $u$, we have
\begin{eqnarray*}
  ({}_xD_L^\mu u, {}_xD_R^\mu u)_{L^2(\Omega)} = \cos(\pi \mu)\|{}_xD_L^\mu u\|^2_{L^2(\Omega)}.
\end{eqnarray*}
\end{Lemma}

\begin{Lemma} (Fractional Poincar{\'e}-Friedrichs, see \cite{e-001}, Theorem 2.10)\label{clp-02-03}
  Let $\mu>0$. For $u\in J_{L,0}^\mu(\Omega)$, we have
\begin{eqnarray*}
  \|u\|_{L^2(\Omega)}\lesssim \|{}_xD_L^\mu u\|_{L^2(\Omega)}.
\end{eqnarray*}
\end{Lemma}

\subsection{Derivation and characterizations of the time-space FE scheme}

Utilizing Lemma \ref{clp-02-01}, referring to \eqref{chp-01-01}, we can derive the variational formulation:
given $\psi_0\in L^2(\Omega)$, $f\in L^2(\Omega,I)$ and $Q_t:=\Omega\times(0,t)$,
to find $u\in \mathcal{H}:=H_0^\gamma(\Omega)\times H^1(I)$ subject to $u(x,0)=\psi_0(x)$ as well as
\begin{equation}\label{chp-02-01}
  \sum_{i=0}^s a_i({}_0^CD_\sigma^{\alpha_i}u,v)_{Q_t}+B^t_\Omega(u,v)=(f,v)_{Q_t},~
  \forall v\in \mathcal{H}^*:=H_0^\gamma(\Omega)\times L^2(I),
\end{equation}
where $({}_0^CD_\sigma^{\alpha_i}u,v)_{Q_t}=\int_0^t({}_0^CD_\sigma^{\alpha_i}u,v)_{L^2(\Omega)}d\sigma$,
$(f,v)_{Q_t}=\int_0^t(f,v)_{L^2(\Omega)}d\sigma$ and
\begin{align*}
  &B^t_\Omega(u,v)=\int_0^t\frac{K_1}{2\cos(\beta\pi)}[({}_xD^\beta_Lu,{}_xD^\beta_Rv)_{L^2(\Omega)}
  +({}_xD^\beta_Ru,{}_xD^\beta_Lv)_{L^2(\Omega)}]d\sigma + \notag \\
  &\qquad\qquad\qquad\qquad\int_0^t\frac{K_2}{2\cos(\gamma\pi)}[({}_xD^\gamma_Lu,{}_xD^\gamma_Rv)_{L^2(\Omega)}
  +({}_xD^\gamma_Ru,{}_xD^\gamma_Lv)_{L^2(\Omega)}]d\sigma.
\end{align*}

To give a description of our time-space FE scheme, we firstly make a temporal mesh $0=t_0<t_1<\cdots<t_N=T$
and a uniform spatial discretization of the interval $\overline{\Omega}$ by points $x_i = a+ih$
($i=0,1,\cdots,M$) with constant spacing $h=(b-a)/M$. We set
\begin{eqnarray*}
  I_j=(t_{j-1},t_j),~\tilde{I}_j=(0,t_j),~\tau_j=t_j-t_{j-1},~j=1,2,\cdots,N;~
  \Omega_h = \{\Omega_l=(x_{l-1},x_l):l=1,2,\cdots,M\}.
\end{eqnarray*}

Next, we introduce the FE spaces in tensor products
\begin{eqnarray*}
  \mathcal{V}_n=\mathcal{V}^\gamma_h(\Omega_h)\times\mathcal{V}_\tau(\tilde{I}_n),~
  \mathcal{V}^*_n=\mathcal{V}^\gamma_h(\Omega_h)\times\mathcal{V}^*_\tau(I_n),
\end{eqnarray*}
where $\mathcal{V}^\gamma_h(\Omega_h)=\{w_h\in H_0^\gamma(\Omega)\cap C(\overline{\Omega}):
w_h(x)|_{\Omega_l}\in\mathcal{P}_1(\Omega_l),~l=1,\cdots,M\}$,
$\mathcal{V}_\tau(\tilde{I}_n) = \{v_\tau\in \mathcal{C}(\overline{\tilde{I}_n}):
v_\tau(0)=1,~v_\tau(t)|_{I_j}\in\mathcal{P}_1(I_j),~j=1,\cdots, n\}$ and
$\mathcal{V}^*_\tau(I_n) = \{v_\tau\in L^2(I_n):v_\tau(t)|_{I_n}\in\mathcal{P}_0(I_n)\}$.
Here $\mathcal{P}_k$ denotes the space of all polynomials of degree $\leq k$.

We now get ready to define the time-space FE numerical scheme of problem \eqref{chp-02-01}: for
$Q_n:=\Omega_h\times I_n$, find $u_{h\tau}\in \mathcal{V}_n$ such that
\begin{eqnarray}\label{chp-02-02}
  \sum_{i=0}^s a_i({}_0^CD_t^{\alpha_i}u_{h\tau},v_{h\tau})_{Q_n}+B^n_\Omega(u_{h\tau},v_{h\tau})=
  (f,v_{h\tau})_{Q_n},~\forall v_{h\tau}\in\mathcal{V}^*_n
\end{eqnarray}
along with $u_{h\tau}(x,0)=\psi_{0,I}(x)$, where $\psi_{0,I}(x)\in\mathcal{V}_n$ satisfies $\psi_{0,I}(x_i)=\psi_0(x_i)$ ($i=0,1,\cdots,M$),
\begin{eqnarray*}
  ({}_0^CD^{\alpha_i}_tu_{h\tau},v_{h\tau})_{Q_n}=\int_{t_{n-1}}^{t_n} ({}_0^CD^{\alpha_i}_t u_{h\tau}, v_{h\tau})_{L^2(\Omega)}dt,~
  (f,v_{h\tau})_{Q_n}=\int_{t_{n-1}}^{t_n}(f,v_{h\tau})_{L^2(\Omega)}dt
\end{eqnarray*}
and
\begin{align*}
  & B^n_\Omega(u_{h\tau},v_{h\tau})
  = \int_{t_{n-1}}^{t_n}\frac{K_1}{2\cos(\beta\pi)}[({}_xD^\beta_Lu_{h\tau},{}_xD^\beta_Rv_{h\tau})_{L^2(\Omega)}
  +({}_xD^\beta_Ru_{h\tau},{}_xD^\beta_Lv_{h\tau})_{L^2(\Omega)}]dt + \notag \\
  &\qquad\qquad\qquad\qquad\int_{t_{n-1}}^{t_n}\frac{K_2}{2\cos(\gamma\pi)}[({}_xD^\gamma_Lu_{h\tau},{}_xD^\gamma_Rv_{h\tau})_{L^2(\Omega)}
  +({}_xD^\gamma_Ru_{h\tau},{}_xD^\gamma_Lv_{h\tau})_{L^2(\Omega)}]dt.
\end{align*}

To go a little further, let
\begin{eqnarray*}
  \mathcal{L}_0(t)=\left \{
    \begin{aligned}
    & \frac{t_1-t}{\tau_1},~t \in I_1 \\
    & 0,~t\in \tilde{I}_n\setminus I_1
    \end{aligned}
  \right.,~
  \mathcal{L}_k(t)=\left \{
    \begin{aligned}
    & \frac{t_{k+1}-t}{\tau_{k+1}},~t \in I_{k+1} \\
    & \frac{t-t_{k-1}}{\tau_k},~t \in I_k \\
    & 0,~t\in \tilde{I}_n\setminus (I_k \cup I_{k+1})
    \end{aligned}
  \right.,~
  \mathcal{L}_n(t)=\left \{
    \begin{aligned}
    & \frac{t-t_{n-1}}{\tau_n},~t \in I_n \\
    & 0,~t\in \tilde{I}_n\setminus I_n
    \end{aligned}
  \right.
\end{eqnarray*}
and
\begin{eqnarray*}
  \tilde{\mathcal{L}}^i_0(t)=\frac{\int_{t_0}^{\hat{t}_1}(t-s)^{-\alpha_i}d\mathcal{L}_0(s)}{\Gamma(1-\alpha_i)},~
  \tilde{\mathcal{L}}^i_k(t)=\frac{\int_{t_{k-1}}^{\hat{t}_{k+1}}(t-s)^{-\alpha_i}d\mathcal{L}_k(s)}{\Gamma(1-\alpha_i)},~
  \tilde{\mathcal{L}}^i_n(t)=\frac{\int_{t_{n-1}}^t(t-s)^{-\alpha_i}d\mathcal{L}_n(s)}{\Gamma(1-\alpha_i)}
\end{eqnarray*}
with $\hat{t}_j=\min(t,t_j)$ for $j=1,\cdots,n$ and $k=1,\cdots,n-1$.

Note that, in view of \eqref{chp-01-02}, the definition of $\mathcal{V}^*_n$ leads to the relation
\begin{eqnarray*}
  \mathcal{V}^*_n=\textsf{span}\{\phi_l(x)\times 1,~l=1,\cdots,M-1\},
\end{eqnarray*}
where $\phi_l(x)$ plays the role of the so-called shape function at point $x_l$. For the monolithic representation
$u_{h\tau}(x,t)=\sum_{k=0}^nu_h^k(x)\mathcal{L}_k(t)$, by direct calculations and taking
\begin{eqnarray*}
  ({}_xD^{\mu}_L\phi_i,{}_xD^{\mu}_R\phi_j)_{L^2(\Omega)} = -({}_xD_L^{2\mu-1}\phi_i,\frac{d\phi_j}{dx})_{L^2(\Omega)},~\mu=\beta,\gamma
\end{eqnarray*}
and the left Riemann-Liouville integral of order $1-2\beta$
\begin{eqnarray*}
  {}_xD_L^{2\beta-1}\phi_i = {}_aJ_x^{1-2\beta}\phi_i := \frac{1}{\Gamma(1-2\beta)}\int_a^x(x-s)^{-2\beta}\phi_i(s)ds
\end{eqnarray*}
into account, one can derive
\begin{align}
  &({}_0^CD^{\alpha_i}_tu_{h\tau},\phi_l\times 1)_{Q_n}
  =\sum_{k=0}^{n-1}(u_h^k,\phi_l)_{L^2(\Omega)}(\tilde{\mathcal{L}}^i_k,1)_{L^2(I_n)}
  +(u_h^n,\phi_l)_{L^2(\Omega)}(\tilde{\mathcal{L}}^i_n,1)_{L^2(I_n)}, \notag \\
  & \int_{t_{n-1}}^{t_n}1\times({}_xD^{\mu}_Lu_{h\tau},{}_xD^{\mu}_R\phi_l)_{L^2(\Omega)}dt
  =\sum_{k=0}^n({}_xD^{\mu}_Lu_h^k,{}_xD^{\mu}_R\phi_l)_{L^2(\Omega)}
  (\mathcal{L}_k,1)_{L^2(I_n)}, \notag \\
  &\int_{t_{n-1}}^{t_n}1\times
  ({}_xD^{\mu}_Ru_{h\tau},{}_xD^{\mu}_L\phi_l)_{L^2(\Omega)}dt
  =\sum_{k=0}^n({}_xD^{\mu}_Ru_h^k,{}_xD^{\mu}_L\phi_l)_{L^2(\Omega)}
  (\mathcal{L}_k,1)_{L^2(I_n)} \notag
\end{align}
and therefore, obtain the desired linear system of equations
\begin{eqnarray}\label{chp-02-03}
  \mathcal{A}^n_{h\tau}\mathcal{U}^n_{h\tau} = \mathcal{F}^n_{h\tau},
\end{eqnarray}
whose coefficient matrix
\begin{eqnarray}\label{chp-02-07}
  \mathcal{A}^n_{h\tau} = \sum_{i=0}^s a_i \frac{\Gamma(3-\alpha_0)\tau_n^{\alpha_0-\alpha_i}}{\Gamma(3-\alpha_i)}M_h
  +K_1\frac{\Gamma(3-\alpha_0)\tau_n^{\alpha_0}}{2}A_h^\beta+K_2\frac{\Gamma(3-\alpha_0)\tau_n^{\alpha_0}}{2}A_h^\gamma,
\end{eqnarray}
right-hand side vector
\begin{align*}
  &\mathcal{F}^n_{h\tau}=\frac{\Gamma(3-\alpha_0)}{\tau_n^{1-\alpha_0}}\bigg{\{}F^n_{h\tau}+
  \Big{[}\sum_{i=0}^s a_i \frac{\tau_n^{1-\alpha_i}}{\Gamma(3-\alpha_i)} M_h-
  K_1\frac{\tau_n}{2} A_h^\beta-K_2\frac{\tau_n}{2} A_h^\gamma\Big{]}
  \mathcal{U}^{n-1}_{h\tau} - \sum_{i=0}^s a_i \times \notag \\
  &\quad \sum_{k=1}^{n-1} \frac{(t_n-t_{k-1})^{2-\alpha_i}-(t_{n-1}-t_{k-1})^{2-\alpha_i}
  -(t_n-t_k)^{2-\alpha_i}+(t_{n-1}-t_k)^{2-\alpha_i}}{\tau_k\Gamma(3-\alpha_i)}
  M_h (\mathcal{U}^k_{h\tau} - \mathcal{U}^{k-1}_{h\tau})\bigg{\}},
\end{align*}
where the vector $F^n_{h\tau}=(f^n_1,f^n_2,\cdots,f^n_{M-1})^T$ with $f^n_l=(f,\phi_l\times 1)_{Q_n}$, the fully FE approximations
\begin{eqnarray*}
  U^k_{h\tau}=(u^k_1,u^k_2,\cdots,u^k_{M-1})^T,~u^0_j=\psi_{0,I}(x_j),~u^k_j=u_h^k(x_j),~j=1,\cdots,M-1,~k=1,\cdots,n,
\end{eqnarray*}
the mass matrix
\begin{eqnarray}\label{chp-02-08}
  M_h = \frac{h}{6}\left(
   \begin{array}{ccccc}
    4 & 1 &   &  &  \\
    1 & 4 & 1 &  &  \\
      & \ddots & \ddots & \ddots &  \\
      &  & 1 & 4 & 1 \\
      &  &   & 1 & 4
   \end{array}
  \right)_{(M-1)\times(M-1)}
\end{eqnarray}
and the stiffness matrices $A_h^\mu = (a^{\mu,h}_{i,j})_{(M-1)\times(M-1)}$ ($\mu=\beta,\gamma$) with each entry of the same form
\begin{align}\label{chp-02-09}
  \left \{
    \begin{aligned}
      & a^{\mu,h}_{i,i} = \frac{h^{1-2\mu}(2^{4-2\mu}-8)}{2\cos(\mu\pi)\Gamma(4-2\mu)},\qquad
      \quad\qquad\qquad\qquad\qquad\qquad\quad ~~~i=1,\cdots,M-1\\
      & a^{\mu,h}_{j,j+1}=a^{\mu,h}_{j+1,j}=\frac{h^{1-2\mu}(3^{3-2\mu}-2^{5-2\mu}+7)}{2\cos(\mu\pi)\Gamma(4-2\mu)},
      \quad\qquad\qquad\quad\quad ~j=1,\cdots,M-2\\
      & a^{\mu,h}_{k,k+l}= a^{\mu,h}_{k+l,k} = \frac{h^{1-2\mu}}{2\cos(\mu\pi)\Gamma(4-2\mu)}[(l+2)^{3-2\mu}\\
      & \quad\quad\quad-4(l+1)^{3-2\mu}+6l^{3-2\mu}-4(l-1)^{3-2\mu}+(l-2)^{3-2\mu}],~k=1,\cdots,M-l-1
    \end{aligned}
  \right..
\end{align}

\begin{Remark}
  The equation \eqref{chp-02-03} follows via multiplying both members of \eqref{chp-02-02} by $\Gamma(3-\alpha_0)\tau_n^{\alpha_0-1}$,
  for the purpose of preventions on severe losses in accuracy and convergence of the time-space FE scheme.
\end{Remark}

An important property on the coefficient matrix $\mathcal{A}^n_{h\tau}$ by \eqref{chp-02-07} is stated in the under-mentioned lemma,
as a natural consequence of the symmetric Toeplitz-like structures of matrices $M_h$ by \eqref{chp-02-08},
$A_h^\beta$ and $A_h^\gamma$ by \eqref{chp-02-09}.

\begin{Lemma}\label{clp-02-05}
  $\mathcal{A}^n_{h\tau}$ is a symmetric Toeplitz matrix.
  Moreover, it is independent of the temporal level $n$ when the time discretization mesh is uniformly spaced.
\end{Lemma}

\begin{Remark}
  Lemma \ref{clp-02-05} implies that (i) a requirement of only $\mathcal{O}(M)$ is used to store $\mathcal{A}^n_{h\tau}$,
  (ii) FFT is the most natural choice for matrix-vector multiplications to
  take advantage of the Toeplitz structure in $\mathcal{O}(M\log M)$ computational complexity.
\end{Remark}

Some important properties of $A_h^\gamma$, the stiffness matrix regarding the diffusion term of \eqref{chp-01-01},
have been already established in \cite{y-001}, which are stated below.

\begin{Lemma} (see \cite{y-001}, Theorem 1) \label{clp-02-04}
  The stiffness matrix $A_h^\gamma$ satisfies
\begin{enumerate}
  \item $a^{\gamma,h}_{i,i} > 0$ for $i=1,\cdots,M-1$;
  \item $a^{\gamma,h}_{i,j} < 0$ for $i\ne j$, $i,j=1,\cdots,M-1$;
  \item $\sum_{j=1}^{M-1} a^{\gamma,h}_{i,j}>0$ for $i=1,\cdots,M-1$.
\end{enumerate}
  As a result, $A_h^\gamma$ is an M-matrix. Moreover, for the particular case when $h\le1/7$, we have
\begin{align}\label{chp-02-04}
  \sum\limits_{j=1}^{M-1} a^{\gamma,h}_{i,j}\ge \left \{
    \begin{aligned}
    & -\frac{h^{1-2\gamma}(4-2^{3-2\gamma})}{2\cos(\gamma\pi)\Gamma(4-2\gamma)},~i=1,M-1 \\
    & -\frac{2^{2\gamma}h(2\gamma-1)}{\cos(\gamma\pi)\Gamma(2-2\gamma)},~~~i=2,\cdots,M-2
    \end{aligned}
  \right..
\end{align}
\end{Lemma}

Now, a few characterizations on $A_h^\beta$, the stiffness matrix from the advection term in \eqref{chp-01-01},
can be also obtained with some important differences.
At this point we denote by $\beta_0$ the unique root of the equation $3^{3-2\beta}-2^{5-2\beta}+7=0$ in the interval (0,1/2).

\begin{Theorem}\label{cth-04-01}
  The stiffness matrix $A_h^\beta$ holds the following properties.
\begin{enumerate}
  \item $a^{\beta,h}_{i,i} > 0$ for $i=1,\cdots,M-1$;
  \item If $\beta>\beta_0$, then $a^{\beta,h}_{i,j} < 0$ for $i\ne j$, $i,j=1,\cdots,M-1$;\\
  Else $a^{\beta,h}_{k,k+1}=a^{\beta,h}_{k+1,k}\ge0$ for $k=1,\cdots,M-2$ and $a^{\beta,h}_{i,j} < 0$ for $|i-j|>1$, $i,j=1,\cdots,M-1$;
  \item $\sum_{j=1}^{M-1} a^{\beta,h}_{i,j}>0$ for $i=1,\cdots,M-1$.
\end{enumerate}
  Thus $A_h^\beta$ is an M-matrix if and only if $\beta\ge\beta_0$. Moreover, for the particular case when $h\le1/7$, we have
\begin{align}\label{chp-02-16}
  \sum\limits_{j=1}^{M-1} a^{\beta,h}_{i,j}\ge \left \{
    \begin{aligned}
    & -\frac{h^{1-2\beta}(4-2^{3-2\beta})}{2\cos(\beta\pi)\Gamma(4-2\beta)},~i=1,M-1 \\
    & -\frac{2^{2\beta}h(2\beta-1)}{\cos(\beta\pi)\Gamma(2-2\beta)},~~~i=2,\cdots,M-2
    \end{aligned}
  \right..
\end{align}
\end{Theorem}

\begin{Proof}
  Property 1 follows immediately since $4-2\beta\in(3,4)$ and $\cos(\beta\pi)>0$.
  It is clear that if $\beta>\beta_0$, then $3^{3-2\beta}-2^{5-2\beta}+7<0$,
  which gives $a^{\beta,h}_{j,j+1}=a^{\beta,h}_{j+1,j}<0$ for $j=1,\cdots,M-2$;
  Otherwise $a^{\beta,h}_{j,j+1}=a^{\beta,h}_{j+1,j}\ge0$ for $j=1,\cdots,M-2$.
  The rest of the property 2 is equivalent to the condition
\begin{eqnarray}\label{chp-02-10}
  f_\beta(l)=(l+2)^{3-2\beta}-4(l+1)^{3-2\beta}+6l^{3-2\beta}-4(l-1)^{3-2\beta}+(l-2)^{3-2\beta} < 0
\end{eqnarray}
  for $2\le l\le M-2$. Hereafter we assume that $\Omega=(0,1)$ without loss of generality. 
  By making use of Taylor's expansion with $x_l=lh$, yields
\begin{align*}
  & f_\beta(l)= h^{1+2\beta}(3-2\beta) (2-2\beta) (1-2\beta) (-2\beta)\Big{[}x_l^{-1-2\beta}+
  \frac{h(-2\beta-1)}{5!} (-4x_{\xi_1}^{-2-2\beta}+4x_{\xi_2}^{-2-2\beta}+\\
  & \quad 32x_{\xi_3}^{-2-2\beta}-32x_{\xi_4}^{-2-2\beta})\Big{]},~~
  x_{\xi_1}\in (x_l,x_{l+1}),~x_{\xi_2}\in (x_{l-1},x_l),~x_{\xi_3}\in (x_l,x_{l+2}),~x_{\xi_4}\in (x_{l-2},x_l).
\end{align*}
  We note that the inequality
\begin{align*}
  & x_l^{-1-2\beta}+\frac{h(2\beta+1)}{5!}(4x_{l+1}^{-2-2\beta}-4x_{l-1}^{-2-2\beta}-32x_{l}^{-2-2\beta}+32x_{l}^{-2-2\beta})=\\
  & \qquad\qquad h^{-1-2\beta}\Big{\{}l^{-1-2\beta}+\frac{4(2\beta+1)}{5!}\Big{[}{(l+1)}^{-2-2\beta}-{(l-1)}^{-2-2\beta}\Big{]} \Big{\}} > 0
\end{align*}
  is clearly true because of the simple observation that the inequality
\begin{eqnarray*}
  (\frac{l}{l+1})^{2+2\beta}-(\frac{l}{l-1})^{2+2\beta} > -\frac{l \times 5!}{4(2\beta+1)}
\end{eqnarray*}
  holds for $2\le l\le M-2$. This gives a derivation of \eqref{chp-02-10}.

  To prove the property 3, it is sufficient to validate
\begin{eqnarray*}
  ({}_xD^{2\beta-1}_L\tilde{\phi},\frac{d\phi_i}{dx})_{L^2(\Omega)}+
  ({}_xD^{2\beta-1}_L\phi_i,\frac{d\tilde{\phi}}{dx})_{L^2(\Omega)}<0,~i=1,\cdots,M-1,
\end{eqnarray*}
  where $\tilde{\phi}:=\sum_{j=1}^{M-1}\phi_j=1-\phi_0-\phi_M$. More precisely,
\begin{eqnarray*}
  \frac{1}{h}\Big{[}({}_xD_L^{2\beta-1}\tilde{\phi},1)_{L^2(\Omega_i)}-({}_xD_L^{2\beta-1}\tilde{\phi},1)_{L^2(\Omega_{i+1})}
  +({}_xD_L^{2\beta-1}\phi_i,1)_{L^2(\Omega_1)}-({}_xD_L^{2\beta-1}\phi_i,1)_{L^2(\Omega_M)}\Big{]} < 0.
\end{eqnarray*}
  Let $\varsigma=3-2\beta\in(2,3)$ for short. For the cases $i=1$ and $i=M-1$, by Taylor's formula, we have to prove
\begin{align*}
  & \hat{f}_\varsigma(i):=\frac{1}{h^2\Gamma(\varsigma+1)}\Big{\{}(4-2^\varsigma)h^\varsigma-\varsigma(\varsigma-1)(\varsigma-2)h^3
  +\frac{3\varsigma(\varsigma-1)(\varsigma-2)(\varsigma-3)}{4!}h^4 \times \\
  & \qquad \Big{[} (1-\xi_1)^{\varsigma-4}-16(1-\xi_2)^{\varsigma-4}+27(1-\xi_3)^{\varsigma-4} \Big{]}\Big{\}} < 0,
  ~~\xi_1 \in (0,h),~\xi_2 \in (0,2h),~\xi_3 \in (0,3h),
\end{align*}
  which is an immediate consequence of the relation
\begin{eqnarray*}
  \frac{(2^\varsigma-4)h^\varsigma}{\varsigma(\varsigma-1)(\varsigma-2)h^3}>\frac{(\varsigma-3)h}{8}[28-16(1-2h)^{\varsigma-4}]-1.
\end{eqnarray*}
  On the other hand, all that is needed is the inequality
\begin{align}\label{chp-02-11}
&\tilde{f}_\varsigma(i):=\frac{1}{h^2\Gamma(\varsigma+1)}\Big{\{}3(ih)^\varsigma-3[(i-1)h]^\varsigma+[(i-2)h]^\varsigma-[(i+1)h]^\varsigma \notag \\
&\qquad \qquad \qquad +3(1-ih)^\varsigma-[1-(i-1)h]^\varsigma-3[1-(i+1)h]^\varsigma+[1-(i+2)h]^\varsigma \Big{\}} < 0.
\end{align}
  In fact, using 6-order Taylor series expansion, it can be easily shown that the inequality
\begin{eqnarray*}
  3(ih)^\varsigma-3[(i-1)h]^\varsigma+[(i-2)h]^\varsigma-[(i+1)h]^\varsigma
  < -h^3\varsigma(\varsigma-1)(\varsigma-2){x_i}^{\varsigma-3}
\end{eqnarray*}
  follows by observing that
\begin{eqnarray*}
  (\frac{i-1}{i})^{6-\varsigma} > \frac{12(5-\varsigma)}{252(5-\varsigma)+i \times 6!},~i=2,\cdots,M-2.
\end{eqnarray*}
  One further can similarly derive
\begin{eqnarray*}
  3(1-ih)^\varsigma-[1-(i-1)h]^\varsigma-3[1-(i+1)h]^\varsigma+[1-(i+2)h]^\varsigma
  < -h^3\varsigma(\varsigma-1)(\varsigma-2){(1-x_i)}^{\varsigma-3}.
\end{eqnarray*}
  Thus, combining the fact that $x_i^{\varsigma-3}+(1-x_i)^{\varsigma-3} \ge 2^{4-\varsigma}$ for any $x_i\in(0,1)$, we obtain
\begin{eqnarray}\label{chp-02-12}
  \tilde{f}_\varsigma(i) < \frac{1}{h^2\Gamma(\varsigma+1)}[-h^3\varsigma(\varsigma-1)(\varsigma-2)][{x_i}^{\varsigma-3}+{(1-x_i)}^{\varsigma-3}]
  \le \frac{h(2-\varsigma)}{\Gamma(\varsigma-1)}2^{4-\varsigma}<0.
\end{eqnarray}
  which completes the proof of \eqref{chp-02-11}.

  Another step is that $A_h^\beta$ is an M-matrix. If $\beta\ge\beta_0$, then $a^{\beta,h}_{j,j+1}=a^{\beta,h}_{j+1,j}\le0$ for $j=1,\cdots,M-2$.
  The converse implication is also true. On this basis, the problem reduces to prove that $(A_h^\beta)^{-1}$ is nonnegative,
  which can be done easily by contradiction, see reference \cite{s-003} for a proof.

  It remains to prove \eqref{chp-02-16}. As the proof of the property 3, for the case when $h\le1/7$, by simple algebraic manipulations, we deduce
\begin{eqnarray*}
  -\varsigma(\varsigma-1)(\varsigma-2)h^3\Big{\{}1-\frac{1}{8}(\varsigma-3)h [(1-\xi_1)^{\varsigma-4}-
  16(1-\xi_2)^{\varsigma-4}+27(1-\xi_3)^{\varsigma-4} ]\Big{\}} < 0,
\end{eqnarray*}
  which leads to $\hat{f}_\varsigma(i) < (4-2^\varsigma)h^{\varsigma-2}/\Gamma(\varsigma+1)$ for $i=1$ and $i=M-1$.
  Furthermore, together with \eqref{chp-02-12}, \eqref{chp-02-16} is proved immediately.
\end{Proof}

For the purpose to ensure that $\mathcal{A}^n_{h\tau}$ is an M-matrix, below are two classes of sufficient conditions.
\begin{description}

  \item[Class 1.] $\beta\ge\beta_0$ and
\begin{eqnarray}\label{chp-02-13}
  \frac{\tau_n^{\alpha_0}}{h^{2\gamma}}>-\frac{4\cos(\gamma\pi)\Gamma(4-2\gamma)}{3K_2(3^{3-2\gamma}
  -2^{5-2\gamma}+7)}\sum_{i=0}^s\frac{a_i}{\Gamma(3-\alpha_i)}.
\end{eqnarray}

  \item[Class 2.] $\beta<\beta_0$ with a suitably small $h$ satisfying
\begin{eqnarray}\label{chp-02-17}
  h^{2(\gamma-\beta)}<-\frac{1}{2} \frac{K_2(3^{3-2\gamma}-2^{5-2\gamma}+7)}{\cos(\gamma\pi)\Gamma(4-2\gamma)}
  \frac{\cos(\beta\pi)\Gamma(4-2\beta)}{K_1(3^{3-2\beta}-2^{5-2\beta}+7)}
\end{eqnarray}
  and \eqref{chp-02-13} concerning a given $\tau_n$.

\end{description}

The upcoming theorem is singled out as an immediate consequence of Lemma \ref{clp-02-04} and Theorem \ref{cth-04-01}.

\begin{Theorem}\label{cth-04-03}
  The following properties are true for the coefficient matrix $\mathcal{A}^n_{h\tau}=(a^{h\tau}_{i,j})_{(M-1)\times(M-1)}$.
\begin{enumerate}

  \item $a^{h\tau}_{i,i} > 0$ for $i=1,\cdots,M-1$;

  \item $a^{h\tau}_{i,j} < 0$ for $|i-j|>1$, $i,j=1,\cdots,M-1$;

  \item It is strictly diagonally dominant.

\end{enumerate}
  In addition, under either Class 1 or Class 2 of sufficient conditions, $\mathcal{A}^n_{h\tau}$ is further an M-matrix.
\end{Theorem}

\begin{Proof}
  Properties 1, 2 and 3 are obvious facts by \eqref{chp-02-07}, \eqref{chp-02-08}, Lemma \ref{clp-02-04} and Theorem \ref{cth-04-01}.
  Evidently, as the proof in Theorem \ref{cth-04-01},
  if we can prove that the minor diagonal elements of $\mathcal{A}^n_{h\tau}$ are all negative
  under either class of sufficient conditions, then the M-matrix property of $\mathcal{A}^n_{h\tau}$ is established immediately.

  In fact, the condition \eqref{chp-02-13} is imposed to make
\begin{eqnarray*}
  \sum_{i=0}^s \frac{a_i\tau_n^{-\alpha_0}}{3\Gamma(3-\alpha_i)}
  +\frac{1}{2}\frac{K_2(3^{3-2\gamma}-2^{5-2\gamma}+7)h^{-2\gamma}}{2\cos(\gamma\pi)\Gamma(4-2\gamma)}<0
\end{eqnarray*}
  valid, which further gives, due to $\tau_n^{\alpha_0-\alpha_i}<1$ ($\alpha_i<\alpha_0$, $i=1,\cdots,s$),
\begin{eqnarray}\label{chp-02-14}
  \sum_{i=0}^s \frac{a_i\tau_n^{-\alpha_i}}{3\Gamma(3-\alpha_i)}
  +\frac{1}{2}\frac{K_2(3^{3-2\gamma}-2^{5-2\gamma}+7)h^{-2\gamma}}{2\cos(\gamma\pi)\Gamma(4-2\gamma)}<0.
\end{eqnarray}

  For the case $\beta\ge\beta_0$, by Theorem \ref{cth-04-01} and \eqref{chp-02-14}, we can conclude that
\begin{eqnarray}\label{chp-02-15}
  \sum_{i=0}^s \frac{a_i\tau_n^{-\alpha_i}}{\Gamma(3-\alpha_i)}\frac{h}{3}
  +\frac{K_1(3^{3-2\beta}-2^{5-2\beta}+7)h^{1-2\beta}}{2\cos(\beta\pi)\Gamma(4-2\beta)}
  +\frac{K_2(3^{3-2\gamma}-2^{5-2\gamma}+7)h^{1-2\gamma}}{2\cos(\gamma\pi)\Gamma(4-2\gamma)}<0.
\end{eqnarray}
  For the opposite $\beta<\beta_0$, we have $(3^{3-2\beta}-2^{5-2\beta}+7)/[2\cos(\beta\pi)\Gamma(4-2\beta)]>0$.
  Observe that the order gap $\gamma-\beta>0$, there always exists a small $h$ subject to \eqref{chp-02-17}.
  Apparently, this situation also yields \eqref{chp-02-15}.
  Furthermore, \eqref{chp-02-15} implies, upon multiplying through the factor $\Gamma(3-\alpha_0)\tau_n^{\alpha_0}/2$,
  that the minor diagonal elements of $\mathcal{A}^n_{h\tau}$ are all negative.
  Hence, results of Theorem \ref{cth-04-03} are obtained.
\end{Proof}

\section{Condition number estimation and an adaptive AMG}

This section is devoted to the derivation of the condition number estimation on $\mathcal{A}^n_{h\tau}$ and the proposal of an appropriate solver.

\subsection{Condition number estimation}

\begin{Theorem}\label{cth-04-02}
  For the matrix $\mathcal{A}^n_{h\tau}$ defined by \eqref{chp-02-07}, we have
\begin{eqnarray}\label{chp-03-01}
  \kappa(\mathcal{A}^n_{h\tau})\lesssim1+\tau_n^{\alpha_0} h^{-2\gamma}.
\end{eqnarray}
\end{Theorem}

\begin{Proof}
  Obviously, there is a spectral equivalence
\begin{eqnarray}\label{chp-03-04}
  \kappa(\mathcal{A}^n_{h\tau})\simeq\kappa\Big{(}\sum_{i=0}^s C_i \tau_n^{\alpha_0-\alpha_i} I + \tilde{C}_1
  \tau_n^{\alpha_0} M_h^{-\frac{1}{2}}A_h^\beta M_h^{-\frac{1}{2}}+\tilde{C}_2 \tau_n^{\alpha_0} M_h^{-\frac{1}{2}}A_h^\gamma M_h^{-\frac{1}{2}}\Big{)},
\end{eqnarray}
  where $C_i=a_i \Gamma(3-\alpha_0) / \Gamma(3-\alpha_i)$ and $\tilde{C}_j=K_j\Gamma(3-\alpha_0)/2$ ($j=1,2$).

  It is worthwhile to point out that the expressions
\begin{eqnarray}\label{chp-03-02}
 \lambda_{\min}(M_h^{-\frac{1}{2}}A_h^\gamma M_h^{-\frac{1}{2}})\gtrsim 1,~
 \lambda_{\max}(M_h^{-\frac{1}{2}}A_h^\gamma M_h^{-\frac{1}{2}})\lesssim h^{-2\gamma}
\end{eqnarray}
  have been obtained by us based on Lemma \ref{clp-02-02}, \ref{clp-02-03} and \ref{clp-02-04}, as well as the Cauchy-Schwarz inequality \cite{y-001}.
  Similarly, according to Lemma \ref{clp-02-02} and \ref{clp-02-03}, we get
\begin{eqnarray*}
  (A_h^\beta \vec{u}_h,\vec{u}_h)=\frac{1}{\cos(\beta\pi)}({}_xD^{\beta}_Lu_h,{}_xD^{\beta}_Ru_h)_{L^2(\Omega)}
  \gtrsim (u_h, u_h)_{L^2(\Omega)}=(M_h\vec{u}_h,\vec{u}_h)\simeq h(\vec{u}_h,\vec{u}_h)
\end{eqnarray*}
  with $u_h:=(\phi_1,\cdots,\phi_{M-1})\vec{u}_h$ for an arbitrary vector $\vec{u}_h:=(u^h_1,\cdots,u^h_{M-1})^T\in \mathbb{R}^{M-1}$.
  Moreover, on account of Theorem \ref{cth-04-01}, we arrive at, for $\beta_0\le\beta<1/2$,
\begin{eqnarray*}
 (A_h^\beta \vec{u}_h, \vec{u}_h) \le \sum_{i=1}^{M-1} a^{\beta,h}_{i,i} (u^h_i)^2 -
 \sum_{i=1}^{M-1} \sum_{j\ne i} a^{\beta,h}_{i,j} \frac{(u^h_i)^2 + (u^h_j)^2}{2}
 \le 2 a^{\beta,h}_{1,1} (\vec{u}_h, \vec{u}_h),
\end{eqnarray*}
  and, for $0<\beta<\beta_0$,
\begin{eqnarray*}
 (A_h^\beta \vec{u}_h, \vec{u}_h) & \le & \sum_{i=1}^{M-1} a^{\beta,h}_{i,i} (u^h_i)^2 +
 \sum_{i=1}^{M-1} \sum_{|j-i|=1} a^{\beta,h}_{i,j} \frac{(u^h_i)^2 + (u^h_j)^2}{2} -
 \sum_{i=1}^{M-1} \sum_{|j-i|>1} a^{\beta,h}_{i,j} \frac{(u^h_i)^2 + (u^h_j)^2}{2} \\
 & = & \sum_{i=1}^{M-1} (u^h_i)^2 \Big{[}a^{\beta,h}_{i,i} + a^{\beta,h}_{i,i+1} +
 a^{\beta,h}_{i,i-1}- \sum_{j=i+2}^{M-1} a^{\beta,h}_{i,j} - \sum_{j=1}^{i-2} a^{\beta,h}_{i,j}\Big{]}
 \le 2 (a^{\beta,h}_{1,1}+2a^{\beta,h}_{1,2}) (\vec{u}_h, \vec{u}_h),
\end{eqnarray*}
  which both indicate that $(A_h^\beta \vec{u}_h,\vec{u}_h)\lesssim h^{1-2\beta}(\vec{u}_h,\vec{u}_h)$,
  regardless of what sort of relationship between $\beta$ and $\beta_0$.
  Taking $\vec{v}_h = M_h^{\frac{1}{2}}\vec{u}_h$, yields
\begin{eqnarray*}
  (\vec{v}_h,\vec{v}_h)\lesssim (M_h^{-\frac{1}{2}}A_h^\beta M_h^{-\frac{1}{2}}\vec{v}_h,\vec{v}_h)
  \lesssim h^{-2\beta}(\vec{v}_h,\vec{v}_h).
\end{eqnarray*}
  As the vector $\vec{v}_h$ is arbitrary, the above expression can be rewritten as
\begin{eqnarray}\label{chp-03-03}
 \lambda_{\min}(M_h^{-\frac{1}{2}}A_h^\beta M_h^{-\frac{1}{2}})\gtrsim 1,~
 \lambda_{\max}(M_h^{-\frac{1}{2}}A_h^\beta M_h^{-\frac{1}{2}})\lesssim h^{-2\beta}.
\end{eqnarray}
  Combining \eqref{chp-03-04}, \eqref{chp-03-02}, \eqref{chp-03-03} and noting that $C_0=a_0>0$, $K_1>0$ and $K_2>0$, give rise to
\begin{eqnarray*}
  \kappa(\mathcal{A}^n_{h\tau})\lesssim\frac{\sum_{i=0}^s C_i \tau_n^{\alpha_0-\alpha_i}+
  \tilde{C}_1\tau_n^{\alpha_0}h^{-2\beta}+\tilde{C}_2 \tau_n^{\alpha_0}h^{-2\gamma}}{\sum_{i=0}^s C_i \tau_n^{\alpha_0-\alpha_i}
  +\tilde{C}_1\tau_n^{\alpha_0}+\tilde{C}_2 \tau_n^{\alpha_0}}
  \lesssim 1+\tau_n^{\alpha_0} h^{-2\gamma},
\end{eqnarray*}
  which completes the proof.
\end{Proof}

\begin{Remark}
  The relation \eqref{chp-03-01} is obviously a generalization of the estimation on resulting matrices
  from traditional parabolic differential equations, single-term time FDEs and FADEs.
\end{Remark}

An immediate corollary of the above theorem follows with great practical interest.

\begin{Corollary}\label{cor-03-01}
  Let the current time step size $\tau_n=\mathcal{O}(h^\varrho)$ and satisfy $\varrho\alpha_0\ge2\gamma$. Then we have
\begin{eqnarray}\label{chp-03-05}
  \kappa(\mathcal{A}^n_{h\tau})=\mathcal{O}(1).
\end{eqnarray}
\end{Corollary}

\subsection{Classical AMG with convergence analysis and its adaptive variant}

Nowadays, classical AMG is quite mature and capable of various ill-conditioned Toeplitz systems.
Most of the existing AMG software packages (e.g. FASP \cite{x-002}, BoomerAMG \cite{h-001}
and AmgX \cite{m-005}) are built on it. It has Setup and Solve phases.
The former phase builds all the ingredients required by a hierarchy of grids, the finest to the coarsest,
while the latter phase runs V-cycle, W-cycle or F-cycle until the desired convergence is achieved,
where the smoother and the coarse-grid correction are crucial components.
It should be emphasized that damped-Jacobi becomes a favorable smoother for the system \eqref{chp-02-03},
since it can be executed by FFT to retain the $\mathcal{O}(M\log M)$ complexity.

Next we turn to the theoretical analysis when $\mathcal{A}^n_{h\tau}$ is a strictly diagonally dominant M-matrix.
We here only assess two-level V(0,1)-cycle, viz., 0 pre-smoothing but 1 post-smoothing step is performed per V-cycle.
It is advantageous to reorder the system \eqref{chp-02-03} and the coarse-to-fine interpolation $\mathcal{P}$
in reference to a given C/F splitting
\begin{eqnarray*}
\mathcal{A}^n_{h\tau}\mathcal{U}^n_{h\tau}=
\left(
\begin{array}{cc}
\mathcal{A}_{FF} & \mathcal{A}_{FC} \\
\mathcal{A}_{CF} & \mathcal{A}_{CC}
\end{array}
\right)
\left(
\begin{array}{c}
\mathcal{U}_F \\
\mathcal{U}_C
\end{array}
\right)
=\left(
\begin{array}{c}
\mathcal{F}_F \\
\mathcal{F}_C
\end{array}
\right) = \mathcal{F}^n_{h\tau},~
\mathcal{P}=\left(
\begin{array}{l}
\mathcal{P}_{FC} \\
\mathcal{I}_{CC}
\end{array}\right)
\end{eqnarray*}
with $\mathcal{I}_{CC}$ as the identity. Another basic tools are specific norms of any vector $v=(v_F,v_C)^T\in \mathbb{R}^{M-1}$
\begin{eqnarray*}
  \|v_F\|_{0,F}=(\textsf{diag}(\mathcal{A}_{FF})v_F,v_F)^{\frac{1}{2}},~
  \|v\|_1=(\mathcal{A}^n_{h\tau}v,v)^{\frac{1}{2}},~
  \|v\|_2=(\textsf{diag}^{-1}(\mathcal{A}^n_{h\tau})\mathcal{A}^n_{h\tau}v,\mathcal{A}^n_{h\tau}v)^{\frac{1}{2}}.
\end{eqnarray*}

It can easily be derived that the two-level iteration matrix to be considered is
\begin{eqnarray*}
  \mathcal{M}_{h,H}=\mathcal{S}_h[\mathcal{I}-\mathcal{P}(\mathcal{P}^T \mathcal{A}^n_{h\tau} \mathcal{P})^{-1}\mathcal{P}^T\mathcal{A}^n_{h\tau}],
\end{eqnarray*}
where $\mathcal{S}_h$ is the post-smoothing iteration matrix.
Using the two-level convergence theory in \cite{t-001} and Theorem \ref{cth-04-03},
the following theorem states an uniform upper bound for $\mathcal{M}_{h,H}$.

\begin{Theorem}\label{cth-03-01}
  Let $\mathcal{S}_h$ be damped-Jacobi or Gauss-Seidel relaxation and $\mathcal{P}$ be the direct interpolation.
  For a given C/F splitting, under either Class 1 or Class 2 of sufficient conditions, there exist constants
  $\sigma_2\ge1>\sigma_1>0$ independent of step sizes $\tau_n$ and $h$, such that
\begin{eqnarray}\label{chp-03-06}
  \|\mathcal{M}_{h,H}\|_1\leq \sqrt{1-\frac{\sigma_1}{\sigma_2}}.
\end{eqnarray}
\end{Theorem}

\begin{Proof}
  By Theorem A.4.1 and A.4.2 in \cite{t-001}, \eqref{chp-03-06} is valid if we provide that $\mathcal{S}_h$ satisfies the smoothing property
\begin{eqnarray}\label{chp-03-07}
  \|\mathcal{S}_he_h\|^2_{1}\leq \|e_h\|^2_1-\sigma_1\|e_h\|^2_2
\end{eqnarray}
  and $\mathcal{P}$ meets the accuracy property
\begin{eqnarray}\label{chp-03-08}
  \|e_F-\mathcal{P}_{FC}e_C\|^2_{0,F} \le \sigma_2\|e_h\|^2_1
\end{eqnarray}
  for all $e_h=(e_F^T,e_C^T)^T\in \mathbb{R}^{M-1}$.

We start with \eqref{chp-03-07}. According to Theorem A.3.1 and A.3.2 in \cite{t-001},
it is valid for damped-Jacobi relaxation (parameter $0<\omega<2/\eta$) with
$\sigma_1=\omega(2-\omega\eta)$,
and for Gauss-Seidel relaxation with
$\sigma_1=1/[(1+\gamma_-)(1+\gamma_+)]$,
where
\begin{eqnarray*}
  \eta \ge \rho(\textsf{diag}^{-1}(\mathcal{A}^n_{h\tau})\mathcal{A}^n_{h\tau}),~
  \gamma_-=\max_i\Big{\{}\frac{1}{w_ia^{h\tau}_{i,i}}\sum_{j<i}w_j|a^{h\tau}_{i,j}|\Big{\}},~
  \gamma_+=\max_i\Big{\{}\frac{1}{w_ia^{h\tau}_{i,i}}\sum_{j>i}w_j|a^{h\tau}_{i,j}|\Big{\}},
\end{eqnarray*}
$\vec{w}=(w_1,\cdots,w_{M-1})^T$ is an arbitrary vector with $\vec{w}$ and $\mathcal{A}^n_{h\tau}\vec{w}$ both being positive.
Recall Theorem \ref{cth-04-03}, one immediately obtains $0<\gamma_-<1$, $0<\gamma_+<1$ and
\begin{eqnarray*}
  \rho(\textsf{diag}^{-1}(\mathcal{A}^n_{h\tau})\mathcal{A}^n_{h\tau})\le|\textsf{diag}^{-1}(\mathcal{A}^n_{h\tau})\mathcal{A}^n_{h\tau}|_{\vec{w}}
  =\max_i\Big{\{}\frac{1}{w_i}\sum_jw_j\frac{|a^{h\tau}_{i,j}|}{a^{h\tau}_{i,i}}\Big{\}}<2.
\end{eqnarray*}
These indicate $\sigma_1\le1/\eta<1$ for damped-Jacobi and $\sigma_1\in(1/4,1)$ for Gauss-Seidel,
both independent of $e_h$, $\tau_n$ and $h$. It is worth noting that Jacobi relaxation is always available,
since there exists a small distance $\epsilon$ so that $\eta=2-\epsilon\ge\rho(\textsf{diag}^{-1}(\mathcal{A}^n_{h\tau})\mathcal{A}^n_{h\tau})$.

It remains to prove the second part \eqref{chp-03-08}. Through Theorem A.4.3 in \cite{t-001}, $\mathcal{P}_{FC}$ holds \eqref{chp-03-08} with
\begin{eqnarray*}
  \sigma_2 \ge \max_{i\in F} \Big{\{}\frac{\sum_{j\in N_i}a^{h\tau}_{i,j}}{\sum_{j\in P_i}a^{h\tau}_{i,j}}\Big{\}},
\end{eqnarray*}
where the neighborhood $N_i=\{j\ne i:a^{h\tau}_{i,j}\ne0\}$, $P_i$ is the set of interpolatory variables at F-point $i$.
Notice the fact that Ruge-St{\"u}ben coarsening strategy retains $i-1\in P_i$ or $i+1\in P_i$ using properties of $\mathcal{A}^n_{h\tau}$.
It implies that $\sigma_2\ge1$ independent of $e_h$ due to $a^{h\tau}_{i,j}<0$ for $j\in N_i$ from Theorem \ref{cth-04-03}, and
\begin{eqnarray}\label{chp-03-09}
  \frac{\sum_{j\in N_i}a^{h\tau}_{i,j}}{\sum_{j\in P_i}a^{h\tau}_{i,j}}<-\frac{a^{h\tau}_{i,i}}{a^{h\tau}_{i,i-1}}=-\frac{a^{h\tau}_{i,i}}{a^{h\tau}_{i,i+1}}
  =-\frac{4C_1+3C_2\tau_n^{\alpha_0}h^{-2\gamma}}{C_1+3C_3\tau_n^{\alpha_0}h^{-2\gamma}},
\end{eqnarray}
where
\begin{eqnarray*}
  C_1=\sum_{i=0}^s a_i \frac{\tau_n^{\alpha_0-\alpha_i}}{\Gamma(3-\alpha_i)},~
  C_2=K_1\frac{h^{2(\gamma-\beta)}(2^{4-2\beta}-8)}{2\cos(\beta\pi)\Gamma(4-2\beta)}
  +K_2\frac{2^{4-2\gamma}-8}{2\cos(\gamma\pi)\Gamma(4-2\gamma)}
\end{eqnarray*}
and
\begin{eqnarray*}
  C_3=K_1\frac{h^{2(\gamma-\beta)}(3^{3-2\beta}-2^{5-2\beta}+7)}{2\cos(\beta\pi)\Gamma(4-2\beta)}
  +K_2\frac{3^{3-2\gamma}-2^{5-2\gamma}+7}{2\cos(\gamma\pi)\Gamma(4-2\gamma)}.
\end{eqnarray*}
Plugging \eqref{chp-02-13}, along with \eqref{chp-02-17} when $\beta<\beta_0$,
\eqref{chp-03-09} suggests that $\sigma_2$ isn't tied to $\tau_n$ and $h$.
Therefore \eqref{chp-03-06} is proved.
\end{Proof}

\begin{Remark}
  The upper bound \eqref{chp-03-06} tells us that a much faster convergence of the two-level V(0,1)-cycle will be achieved if
  $\sigma_1/\sigma_2$ can be made quite closer to 1.
\end{Remark}

Table \ref{ctp-03-01} is offered in support of Theorem \ref{cth-03-01}.
It is easy to see that the uniform convergence of the two-level V(0,1)-cycle is achieved.
Nevertheless, the strength-of-connection tolerance $\theta$, used to interpret concepts of strong influence and dependence,
has a significantly negative impact on its actual number of iterations.

\begin{table}[htbp]
\footnotesize
\centering\caption{Effect of $\theta$ on the two-level V(0,1)-cycle with $10^{-8}$ as the tolerance for stopping.}\label{ctp-03-01}\vskip 0.1cm
\begin{tabular}{||c|c|c|c|c|c|c|c|c|c|c|}\hline
\multirow{2}{*}{$M$}&\multicolumn{3}{c|}{$\alpha_0=0.9$, $\alpha_1=0.4$, $\beta=0.3$, $\gamma=0.8$}&
\multicolumn{3}{c||}{$\alpha_0=0.7$, $\alpha_1=0.5$, $\beta=0.15$, $\gamma=0.95$} \\ \cline{2-7}
~& $\theta=0.0001$ & $\theta=0.005$ & $\theta=0.25$ & $\theta=0.0001$ & $\theta=0.005$ & \multicolumn{1}{c||}{$\theta=0.25$} \\\hline

 512 & 292 & 36 & 13 & 191 & 35 & \multicolumn{1}{c||}{15} \\ \hline
1024 & 308 & 38 & 13 & 198 & 35 & \multicolumn{1}{c||}{15} \\ \hline
2048 & 314 & 36 & 13 & 202 & 35 & \multicolumn{1}{c||}{15} \\ \hline

\end{tabular}
\end{table}

It is well to be reminded that multigrid V(1,1)-cycle in CF-relaxation is of a larger practical value in applications. In this situation,
the impact becomes rather serious. Table \ref{ctp-03-02} shows the results of $\alpha_0=0.7$, $\alpha_1=0.5$, $\beta=0.15$, $\gamma=0.95$.

\begin{table}[htbp]
\footnotesize
\centering\caption{Effect of $\theta$ on multigrid V(1,1)-cycle
in CF-relaxation with $10^{-12}$ as the tolerance for stopping.}\label{ctp-03-02}\vskip 0.1cm
\begin{tabular}{||c|c|c|c|c|c|c|c|c|c|c|c|c|c|c|c|c|c|c|}\hline
\multirow{2}{*}{$M$}& \multicolumn{3}{c|}{$\theta=0.0001$} & \multicolumn{3}{c|}{$\theta=0.03529$}
& \multicolumn{3}{c||}{$\theta=0.03533$} \\ \cline{2-10}
~& Levels & Iterations & CPU & Levels & Iterations & CPU & Levels & Iterations & \multicolumn{1}{c||}{CPU} \\\hline

 512 & 3 & 224 & 6.08E-1 & 6 & 21 & 1.38E-1 & 7 & 5 & \multicolumn{1}{c||}{3.80E-2} \\ \hline
1024 & 3 & 226 & 2.53E0  & 7 & 23 & 5.42E-1 & 8 & 5 & \multicolumn{1}{c||}{1.55E-1} \\ \hline
2048 & 4 & 227 & 1.02E1  & 8 & 25 & 2.23E0  & 9 & 5 & \multicolumn{1}{c||}{6.17E-1} \\ \hline

\end{tabular}
\end{table}

As indicated, there are two drawbacks to multigrid V(1,1)-cycle
in CF-relaxation. One lies the good choice of $\theta$. We have discussed the most reliable
$\theta=a^{h\tau}_{1,3}/a^{h\tau}_{1,2}+\epsilon_0$ for time-space Caputo-Riesz FDEs \cite{y-001},
where $\epsilon_0$ is some small number. However, it may cause the method lack of
convergence optimality because of positive minor diagonal elements of $A^\beta_h$ when $\beta<\beta_0$,
e.g. in Table \ref{ctp-03-02}, $a^{h\tau}_{1,3}/a^{h\tau}_{1,2}\thickapprox0.035285$ for the case $M=512$,
but 0.03529 isn't the optimal value of $\theta$. A natural cure technique for this is to employ
\begin{eqnarray}\label{chp-03-10}
  \theta^{(k)}=\frac{a^{(k)}_{1,3}}{a^{(k)}_{1,2}}+\epsilon_0
\end{eqnarray}
for the $k$-th coarse-level Galerkin matrix $\mathcal{A}_k=(a^{(k)}_{i,j})$ with $\mathcal{A}_0=\mathcal{A}^n_{h\tau}$.
Indeed, $\max_k(\theta^{(k)})\thickapprox0.03533$ to obtain the optimal convergence.
The other is the required cost of $\mathcal{O}(M^2)$ at each time step.
This drawback hampers the acceptance of the method in large linear systems.
The primary reason is that it is impossible to make
$\mathcal{A}_1=\mathcal{P}^T\mathcal{A}_0\mathcal{P}$ Toeplitz-like by direct interpolation in a straightforward way.
This requests for a modification to $\mathcal{P}$. The following lemma provides a proper manipulation.

\begin{Lemma}\label{clp-03-01}
  If $\mathcal{P}_{FC}$ admits 1/2 as all nonzero entries based on the choice \eqref{chp-03-10},
  then $\mathcal{A}_1$ must be a symmetric Toeplitz matrix.
\end{Lemma}

\begin{Proof}
  Obviously, $\mathcal{A}_1^T = (\mathcal{P}^T\mathcal{A}_0\mathcal{P})^T = \mathcal{P}^T\mathcal{A}_0\mathcal{P} = \mathcal{A}_1$, and
\begin{eqnarray*}
  a^{(1)}_{i,i+l}=a^{(1)}_{i+l,i}=\frac{1}{4}a^{h\tau}_{3,2l+1}+a^{h\tau}_{2,2l+1}+
  \frac{3}{2}a^{h\tau}_{1,2l+1}+a^{h\tau}_{1,2l+2}+\frac{1}{4}a^{h\tau}_{1,2l+3},~i=1,\cdots,M_1-l,
\end{eqnarray*}
  where $M_1$ is the number of columns of $\mathcal{P}_{FC}$. This completes the proof.
\end{Proof}

\begin{Remark}
  It can easily be seen that just $\mathcal{O}(M_1)$ operations are required to generate $\mathcal{A}_1$
  (only $M_1$ different entries), whose cost is negligible relative to that of direct interpolation.
\end{Remark}

More interestingly, the two-level convergence for the modified interpolation is also uniform as a result of the following theorem.

\begin{Theorem}
  Let $h\le1/7$ and $\mathcal{P}_{FC}$ satisfy Lemma \ref{clp-03-01}. For a given C/F splitting,
  under either Class 1 or Class 2 of sufficient conditions, there exists a constant
  $\sigma_2\ge1$ independent of step sizes $\tau_n$ and $h$,
  such that the accuracy property \eqref{chp-03-08} holds for all $e_h=(e_F^T,e_C^T)^T=(e_1,e_2,\cdots,e_{M-1})^T\in \mathbb{R}^{M-1}$.
\end{Theorem}

\begin{Proof}
To estimate the right part, we note
\begin{eqnarray}\label{chp-03-11}
  \|e_h\|^2_1=(\mathcal{A}^n_{h\tau}e_h,e_h)\ge\sum_{i\in F}\Big{[}\sum_{k\in P_i}(-a^{h\tau}_{i,k})(e_i-e_k)^2
  +\sum_ja^{h\tau}_{i,j}e_i^2\Big{]}.
\end{eqnarray}
On the other hand, we can estimate
\begin{eqnarray}\label{chp-03-12}
  \|e_F-\mathcal{P}_{FC}e_C\|^2_{0,F}=\sum_{i\in F}a^{h\tau}_{i,i}\Big{(}e_i-\sum_{k\in P_i}\frac{1}{2}e_k\Big{)}^2
  \le\sum_{i\in F}a^{h\tau}_{i,i}\sum_{k\in P_i}\frac{1}{2}(e_i-e_k)^2+a^{h\tau}_{1,1}(e_1^2+e_{M-1}^2),
\end{eqnarray}
employing Schwarz inequality and the fact that Ruge-St{\"u}ben coarsening strategy guarantees $i-1$ or $i+1$ inside of set $P_i$ at F-point $i$.
The previous estimations \eqref{chp-03-11} and \eqref{chp-03-12} imply \eqref{chp-03-08}, if relations
$\sigma_2 (-a^{h\tau}_{i,k}) \ge a^{h\tau}_{i,i}/2$ as well as $\sigma_2\sum_{j=1}^{M-1}a^{h\tau}_{1,j}\ge a^{h\tau}_{1,1}$
hold simultaneously for $k\in P_i$ and $i\in F$, namely
\begin{eqnarray}\label{chp-03-09}
  \sigma_2\ge\max_{i\in F}\Big{\{}-\frac{a^{h\tau}_{i,i}}{2a^{h\tau}_{i,i-1}},
  \frac{a^{h\tau}_{1,1}}{\sum_{j=1}^{M-1}a^{h\tau}_{1,j}}\Big{\}}.
\end{eqnarray}
It indicates that $\sigma_2\ge1$ independent of $e_h$. Furthermore, observing \eqref{chp-02-07}, \eqref{chp-02-08},
\eqref{chp-02-04} and \eqref{chp-02-16}, yield
\begin{eqnarray*}
  \frac{a^{h\tau}_{1,1}}{\sum_{j=1}^{M-1}a^{h\tau}_{1,j}}<2,~
  -\frac{a^{h\tau}_{i,i}}{2a^{h\tau}_{i,i-1}}=-\frac{4C_1+3(K_1C^\beta_2h^{2(\gamma-\beta)}+K_2C^\gamma_2)\tau_n^{\alpha_0}
  h^{-2\gamma}}{2C_1+6(K_1C^\beta_3h^{2(\gamma-\beta)}+K_2C^\gamma_3)\tau_n^{\alpha_0}h^{-2\gamma}},
\end{eqnarray*}
where
\begin{eqnarray*}
  C_1=\sum_{i=0}^s a_i \frac{\tau_n^{\alpha_0-\alpha_i}}{\Gamma(3-\alpha_i)},~
  C^\varrho_2=\frac{2^{4-2\varrho}-8}{2\cos(\varrho\pi)\Gamma(4-2\varrho)},~
  C^\varrho_3=\frac{3^{3-2\varrho}-2^{5-2\varrho}+7}{2\cos(\varrho\pi)\Gamma(4-2\varrho)}.
\end{eqnarray*}
Plugging \eqref{chp-02-13}, along with \eqref{chp-02-17} when $\beta<\beta_0$,
\eqref{chp-03-09} suggests that $\sigma_2$ isn't tied to $\tau_n$ or $h$ at all.
This establishes the result.
\end{Proof}

It should be stressed that these previous procedures can be performed similarly with multigrid
to render all $\mathcal{A}_k$ ($k\ge2$) Toeplitz-like.
Two practical benefits involve computational cost of $\mathcal{O}(M\log M)$ and matrix-free storage of $\mathcal{O}(M)$.
Additionally, it is easy to verify that $\mathcal{A}_k$ ($k\ge1$) could also be an M-matrix for small $k$
constrained by \eqref{chp-02-13} with $h$ being replaced by $2^kh$, together with \eqref{chp-02-17} if $\beta<\beta_0$.

We conclude this section by developing an adaptive AMG algorithm, which exploits features of FFT,
Toeplitz-like structures of $\mathcal{A}_k$ ($k\ge0$) and at most 2 nonzeros per row in interpolations,
while Corollary \ref{cor-03-01} as our clear distinction.

\begin{Algorithm} \label{alg-03-01}
An adaptive AMG for the system \eqref{chp-02-03}.
\begin{description}

  \item[Step 1.] If $\mathcal{A}^n_{h\tau}$ satisfies the estimation \eqref{chp-03-05},
  then solve \eqref{chp-02-03} by conjugate gradient (CG) algorithm.

  \item[Step 2.] Apply the improved classical AMG to solve \eqref{chp-02-03} until convergence.
  \begin{description}

    \item[2.1] Setup phase.
    \begin{description}

      \item[2.1.1] Set $\epsilon_0=10^{-8}$, \textsf{max\_cdofs}, \textsf{max\_levels}, $j=0$ and $\Omega^{(0)}=\{1,2,\cdots,M-1\}$.

      \item[2.1.2] Compute $\theta^{(j)}$ via \eqref{chp-03-10}, perform $\Omega^{(j)}=C^{(j)}\cup F^{(j)}$,
      and set $\Omega^{(j+1)}=C^{(j)}$.

      \item[2.1.3] Construct the modified $\mathcal{P}^{(j)}$
      and $\mathcal{A}_{j+1}=(\mathcal{P}^{(j)})^T\mathcal{A}_j\mathcal{P}^{(j)}$.

      \item[2.1.4] If $|\Omega^{(j+1)}|\le$\textsf{max\_cdofs} or $j+1=$\textsf{max\_levels}, then Stop; else $j=j+1$, goto 2.1.2.

    \end{description}

    \item[2.2] Solve phase: V(1,1)-cycle.
    Set $\mathcal{F}_0=\mathcal{F}^n_{h\tau}$ and choose an initial guess $\mathcal{U}_0$.
    \begin{description}

      \item[2.2.1] For $k=0,1,\cdots,j$, do:
      \begin{itemize}
        \item Run Jacobi relaxation once to $\mathcal{A}_k\mathcal{U}_k = \mathcal{F}_k$ with $\mathcal{U}_k=0$ ($k>0$).
        \item Compute $\mathcal{F}_{k+1}=(\mathcal{P}^{(k)})^T(\mathcal{F}_k-\mathcal{A}_k\mathcal{U}_k)$.
      \end{itemize}

      \item[2.2.2] Solve $\mathcal{A}_{j+1}\mathcal{U}_{j+1} = \mathcal{F}_{j+1}$ by using Gaussian elimination without pivoting.

      \item[2.2.3] For $k=j,j-1,\cdots,0$, do:
      \begin{itemize}
        \item Update $\mathcal{U}_k=\mathcal{U}_k+\mathcal{P}^{(k)}\mathcal{U}_{k+1}$.
        \item Run Jacobi relaxation once to $\mathcal{A}_k\mathcal{U}_k = \mathcal{F}_k$.
      \end{itemize}

    \end{description}

  \end{description}

\end{description}
\end{Algorithm}

\section{Numerical results}

In this section we will present some numerical experiments to illustrate the convergence of our time-space FE numerical scheme \eqref{chp-02-02},
the correctness of our condition number estimation and the effectiveness of our solver.
All examples are implemented in C on a 64 bit Fedora 18 platform with an -O2 optimization parameter,
double precision arithmetic on Intel Xeon (W5590) with configurations 24.0 GB RAM and 3.33 GHz.

In our implementations, all integrals are calculated by a quadrature formula. In tables below,
columns labeled $\|e\|_0$ represent errors $\|u(\cdot,T)-u_{h\tau}(\cdot,T)\|_{L^2(\Omega)}$,
$\lambda_{\min}$ denote the smallest eigenvalues and $\lambda_{\max}$ as the largest eigenvalues,
$\kappa$ are condition numbers of resulting matrices,
\emph{Its} indicate numbers of iterations until the actual residual is reduced by a factor of $10^{-12}$,
$T_c$ express CPU times including Setup and Solve phases with second as the unit,
starred entries ($\ast$) indicate the solutions fail to converge after 1000 iterations,
while CG, CAMG and $i$CAMG stand for conjugate gradient, classical AMG and the adaptive variant in Algorithm \ref{alg-03-01}, respectively.

\subsection{Convergence test}

\begin{Example}\label{cex-04-01}
  Consider problem \eqref{chp-01-01}-\eqref{chp-01-03} with $s=1$, $a_0=a_1=1$, $K_1=1$, $K_2=2$, $T=0.5$, $\Omega=(0,1)$, $\psi_0(x)=100(x^2-x^3)$ and
\begin{align*}
  & f(x,t) =200(x^2-x^3)\Big{[}\frac{t^{2-\alpha_0}}{\Gamma(3-\alpha_0)}+\frac{t^{2-\alpha_1}}{\Gamma(3-\alpha_1)}\Big{]} +
  \frac{50(t^2+1)}{\cos(\beta\pi)}\times\Big{[}\frac{(1-x)^{1-2\beta}}{\Gamma(2-2\beta)}+\\
  &\qquad\qquad\qquad\frac{2x^{2-2\beta}-4(1-x)^{2-2\beta}}{\Gamma(3-2\beta)} + \frac{6(1-x)^{3-2\beta}-6x^{3-2\beta}}{\Gamma(4-2\beta)}\Big{]}+
  \frac{100(t^2+1)}{\cos(\gamma\pi)}\times\\
  &\qquad\qquad\qquad\qquad\Big{[}\frac{(1-x)^{1-2\gamma}}{\Gamma(2-2\gamma)}
  +\frac{2x^{2-2\gamma}-4(1-x)^{2-2\gamma}}{\Gamma(3-2\gamma)} +
  \frac{6(1-x)^{3-2\gamma}-6x^{3-2\gamma}}{\Gamma(4-2\gamma)}\Big{]}.
\end{align*}
\end{Example}

\begin{figure}[htp]
\begin{minipage}[t]{0.5\linewidth}
 \centering
 \subfigure[The exact solution]{%
 \label{fig-04-01-a}
 \includegraphics[width=2.95in]{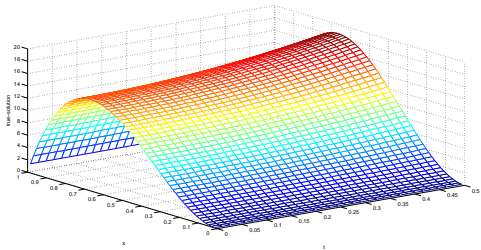}}
\end{minipage}%
\begin{minipage}[t]{0.5\linewidth}
 \centering
 \subfigure[The numerical solution]{%
 \label{fig-04-01-b}
 \includegraphics[width=2.95in]{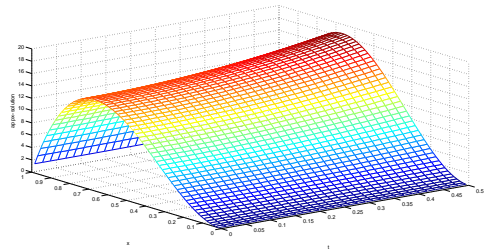}}
\end{minipage}
\caption{Comparison on solutions, $h=\tau=1/64$,
$\alpha_0=0.7$, $\alpha_1=0.4$, $\beta=0.15$, $\gamma=0.95$.}\label{fig-04-01}
\end{figure}

The exact solution is $u(x,t)=100(t^2+1)(x^2-x^3)$. Table \ref{ctp-04-01} and \ref{ctp-04-02} show the results of errors and convergence rates
with typical $\alpha_0$, $\alpha_1$, $\beta$ and $\gamma$ for two specific cases: $h=\tau$ and $h=\sqrt{\tau}$, respectively.
An interpretation is that the time-space FE solution possesses the saturation error order $\mathcal{O}(\tau^2+h^2)$.
Fig. \ref{fig-04-01} illustrates the exact solution with the numerical solution of
$\alpha_0=0.7$, $\alpha_1=0.4$, $\beta=0.15$, $\gamma=0.95$ when $h=\tau=1/64$.

\begin{table}[htbp]
\small
\centering\caption{Error results and convergence rates in spatial direction with $h=\tau$.}\label{ctp-04-01}\vskip 0.1cm
\begin{tabular}{||c|c|c|c|c|c|c|c|c|c|c|c|c|c|c|c|c|c|c|}\hline
\multirow{3}{*}{$M$}&\multicolumn{4}{c|}{$\alpha_0=0.5$, $\alpha_1=0.2$}&
\multicolumn{4}{c||}{$\alpha_0=0.7$, $\alpha_1=0.4$} \\ \cline{2-9}
~&\multicolumn{2}{c|}{$\beta=0.3$, $\gamma=0.8$}&\multicolumn{2}{c|}{$\beta=0.15$, $\gamma=0.95$}&
\multicolumn{2}{c|}{$\beta=0.3$, $\gamma=0.8$}&\multicolumn{2}{c||}{$\beta=0.15$, $\gamma=0.95$} \\ \cline{2-9}
~& $\|e\|_0$ & rate & $\|e\|_0$ & rate & $\|e\|_0$ & rate & $\|e\|_0$ & \multicolumn{1}{c||}{rate} \\\hline

 16 & 6.837E-2 &   -   & 8.357E-2 &   -   & 6.396E-2 &   -   & 8.186E-2 & \multicolumn{1}{c||}{  -  } \\ \hline
 32 & 1.525E-2 & 2.165 & 2.020E-2 & 2.049 & 1.458E-2 & 2.133 & 1.981E-2 & \multicolumn{1}{c||}{2.047} \\ \hline
 64 & 3.484E-3 & 2.130 & 4.878E-3 & 2.050 & 3.383E-3 & 2.108 & 4.811E-3 & \multicolumn{1}{c||}{2.042} \\ \hline
128 & 8.113E-4 & 2.102 & 1.183E-3 & 2.044 & 7.948E-4 & 2.089 & 1.171E-3 & \multicolumn{1}{c||}{2.039} \\ \hline

\end{tabular}
\end{table}

\begin{table}[htbp]
\small
\centering\caption{Error results and convergence rates in spatial direction with $h=\sqrt{\tau}$.}\label{ctp-04-02}\vskip 0.1cm
\begin{tabular}{||c|c|c|c|c|c|c|c|c|c|c|c|c|c|c|c|c|c|c|}\hline
\multirow{3}{*}{$N$}&\multicolumn{4}{c|}{$\alpha_0=0.5$, $\alpha_1=0.2$}
&\multicolumn{4}{c||}{$\alpha_0=0.7$, $\alpha_1=0.4$} \\ \cline{2-9}
~&\multicolumn{2}{c|}{$\beta=0.3$, $\gamma=0.8$}&\multicolumn{2}{c|}{$\beta=0.15$, $\gamma=0.95$}&
\multicolumn{2}{c|}{$\beta=0.3$, $\gamma=0.8$}&\multicolumn{2}{c||}{$\beta=0.15$, $\gamma=0.95$} \\ \cline{2-9}
~& $\|e\|_0$ & rate & $\|e\|_0$ & rate & $\|e\|_0$ & rate & $\|e\|_0$ & \multicolumn{1}{c||}{rate} \\\hline

 32  & 2.600E-1 &   -   & 3.182E-1 &   -   & 2.582E-1 &   -   & 3.164E-1 & \multicolumn{1}{c||}{  -  } \\ \hline
 128 & 5.929E-2 & 1.066 & 7.719E-2 & 1.022 & 5.899E-2 & 1.065 & 7.692E-2 & \multicolumn{1}{c||}{1.020} \\ \hline
 512 & 1.369E-2 & 1.057 & 1.877E-2 & 1.020 & 1.362E-2 & 1.058 & 1.871E-2 & \multicolumn{1}{c||}{1.020} \\ \hline
2048 & 3.194E-3 & 1.050 & 4.569E-3 & 1.019 & 3.177E-3 & 1.050 & 4.554E-3 & \multicolumn{1}{c||}{1.019} \\ \hline

\end{tabular}
\end{table}

\begin{Example}\label{cex-04-02}
  Consider problem \eqref{chp-01-01}-\eqref{chp-01-03} with $s=1$, $a_0=a_1=1$, $T=0.5$, $\Omega=(0,1)$, $\psi_0(x)=100x^2(1-x)^2$ and
\begin{align*}
  & f(x,t) =200x^2(1-x)^2\Big{[}\frac{t^{2-\alpha_0}}{\Gamma(3-\alpha_0)}+\frac{t^{2-\alpha_1}}{\Gamma(3-\alpha_1)}\Big{]} +
  \frac{100K_1(t^2+1)}{\cos(\beta\pi)}\times\Big{[}\frac{x^{2-2\beta}+(1-x)^{2-2\beta}}{\Gamma(3-2\beta)}-\\
  &\qquad\qquad\qquad\frac{6x^{3-2\beta}+6(1-x)^{3-2\beta}}{\Gamma(4-2\beta)} + \frac{12x^{4-2\beta}+12(1-x)^{4-2\beta}}{\Gamma(5-2\beta)}\Big{]}+
  \frac{100K_2(t^2+1)}{\cos(\gamma\pi)}\times\\
  &\qquad\qquad\qquad\qquad\Big{[}\frac{x^{2-2\gamma}+(1-x)^{2-2\gamma}}{\Gamma(3-2\gamma)}
  -\frac{6x^{3-2\gamma}+6(1-x)^{3-2\gamma}}{\Gamma(4-2\gamma)} +
  \frac{12x^{4-2\gamma}+12(1-x)^{4-2\gamma}}{\Gamma(5-2\gamma)}\Big{]}.
\end{align*}
\end{Example}

\begin{figure}[htp]
\begin{minipage}[t]{0.5\linewidth}
 \centering
 \subfigure[The exact solution]{%
 \label{fig-04-02-a}
 \includegraphics[width=2.95in]{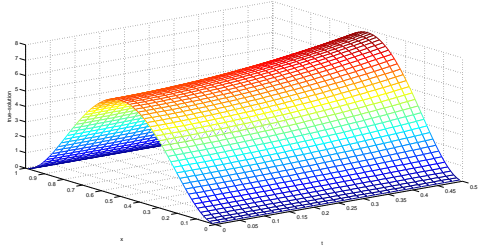}}
\end{minipage}%
\begin{minipage}[t]{0.5\linewidth}
 \centering
 \subfigure[The numerical solution]{%
 \label{fig-04-02-b}
 \includegraphics[width=2.95in]{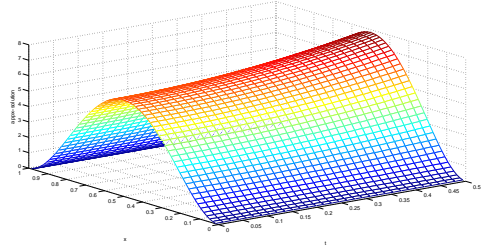}}
\end{minipage}
\caption{Comparison on solutions, $h=\tau=1/64$, $K_1=5$, $K_2=30$,
$\alpha_0=0.7$, $\alpha_1=0.4$, $\beta=0.3$, $\gamma=0.85$.}\label{fig-04-02}
\end{figure}

The exact solution is $u(x,t)=100(t^2+1)x^2(1-x)^2$.
The objective of this example is to measure the possible effects of $K_1$ and $K_2$ on the convergence.
Setting $K_1=5$ and $K_2=30$, $300$, $10^3$ and $10^6$, we can observe from Table \ref{ctp-04-03} and \ref{ctp-04-04} that time-space FE solutions
retain $\mathcal{O}(\tau^2+h^2)$, without any ties to $K_1$ nor $K_2$. Fig. \ref{fig-04-02} charts the exact and numerical solutions distributed over
$[0,1]\times[0,0.5]$ when $K_1=5$, $K_2=30$. This comparison demonstrates that the time-space FE solution confirms the exact solution very well.

\begin{table}[htbp]
\small
\centering\caption{Error results and convergence rates in spatial direction with
$\alpha_0=0.7$, $\alpha_1=0.4$, $\beta=0.3$, $\gamma=0.85$.}\label{ctp-04-03}\vskip 0.1cm
\begin{tabular}{||c|c|c|c|c|c|c|c|c|c|c|c|c|c|c|c|c|c|c|}\hline
\multirow{3}{*}{$M$}&\multicolumn{4}{c|}{$h=\tau$}
&\multirow{3}{*}{$N$}&\multicolumn{4}{c||}{$h=\sqrt{\tau}$} \\ \cline{2-5} \cline{7-10}
~&\multicolumn{2}{c|}{$K_1=5$, $K_2=30$}&\multicolumn{2}{c|}{$K_1=5$, $K_2=300$}&~&
\multicolumn{2}{c|}{$K_1=5$, $K_2=30$}&\multicolumn{2}{c||}{$K_1=5$, $K_2=300$} \\ \cline{2-5} \cline{7-10}
~& $\|e\|_0$ & rate & $\|e\|_0$ & rate & ~& $\|e\|_0$ & rate & $\|e\|_0$ & \multicolumn{1}{c||}{rate} \\\hline

 16 & 3.455E-2 &   -   & 3.607E-2 &   -   &  32  & 1.305E-1 &   -   & 1.368E-1 & \multicolumn{1}{c||}{  -  } \\ \hline
 32 & 8.466E-3 & 2.029 & 8.774E-3 & 2.040 &  128 & 3.065E-2 & 1.045 & 3.302E-2 & \multicolumn{1}{c||}{1.026} \\ \hline
 64 & 1.987E-3 & 2.090 & 2.121E-3 & 2.049 &  512 & 7.362E-3 & 1.029 & 7.509E-3 & \multicolumn{1}{c||}{1.068} \\ \hline
128 & 4.509E-4 & 2.140 & 5.228E-4 & 2.020 & 2048 & 1.774E-3 & 1.027 & 1.792E-3 & \multicolumn{1}{c||}{1.034} \\ \hline

\end{tabular}
\end{table}

\begin{table}[htbp]
\small
\centering\caption{Error results and convergence rates in spatial direction with
$\alpha_0=0.8$, $\alpha_1=0.3$, $\beta=0.2$, $\gamma=0.75$.}\label{ctp-04-04}\vskip 0.1cm
\begin{tabular}{||c|c|c|c|c|c|c|c|c|c|c|c|c|c|c|c|c|c|c|}\hline
\multirow{3}{*}{$M$}&\multicolumn{4}{c|}{$h=\tau$}
&\multirow{3}{*}{$N$}&\multicolumn{4}{c||}{$h=\sqrt{\tau}$} \\ \cline{2-5} \cline{7-10}
~&\multicolumn{2}{c|}{$K_1=5$, $K_2=10^3$}&\multicolumn{2}{c|}{$K_1=5$, $K_2=10^6$}&~&
\multicolumn{2}{c|}{$K_1=5$, $K_2=10^3$}&\multicolumn{2}{c||}{$K_1=5$, $K_2=10^6$} \\ \cline{2-5} \cline{7-10}
~& $\|e\|_0$ & rate & $\|e\|_0$ & rate & ~& $\|e\|_0$ & rate & $\|e\|_0$ & \multicolumn{1}{c||}{rate} \\\hline

 16 & 3.544E-2 &   -   & 3.571E-2 &   -   &  32  & 1.350E-01 &   -   & 1.418E-01 & \multicolumn{1}{c||}{  -  } \\ \hline
 32 & 8.635E-3 & 2.037 & 8.844E-3 & 2.013 &  128 & 3.221E-02 & 1.034 & 3.567E-02 & \multicolumn{1}{c||}{0.996} \\ \hline
 64 & 2.065E-3 & 2.064 & 2.187E-3 & 2.016 &  512 & 6.905E-03 & 1.111 & 8.720E-03 & \multicolumn{1}{c||}{1.016} \\ \hline
128 & 4.916E-4 & 2.071 & 5.415E-4 & 2.014 & 2048 & 1.589E-03 & 1.060 & 2.047E-03 & \multicolumn{1}{c||}{1.045} \\ \hline

\end{tabular}
\end{table}

\subsection{Condition number test}

\begin{Example}\label{cex-04-03}
  We use the same Example \ref{cex-04-01}.
\end{Example}

In this test, we examine the effectiveness of our estimation \eqref{chp-03-01} in three situations: $\tau=h$, $\tau=h^2$ and $\tau=1/64$.
Results are summarized in Table \ref{ctp-04-05}, \ref{ctp-04-06} and \ref{ctp-04-07}.
We observe that they all coincide uniformly with Theorem \ref{cth-04-02}.

\begin{table}[htbp]
\small
\centering\caption{The smallest and largest eigenvalues, condition numbers and ratios with $\tau=h$.}\label{ctp-04-05}\vskip 0.1cm
\begin{tabular}{||c|c|c|c|c|c|c|c|c|c|}\hline
\multirow{2}{*}{$M$}&\multicolumn{4}{c|}{$\alpha_0=0.9$, $\alpha_1=0.4$, $\beta=0.3$, $\gamma=0.8$}
&\multicolumn{4}{c||}{$\alpha_0=0.7$, $\alpha_1=0.5$, $\beta=0.15$, $\gamma=0.95$} \\ \cline{2-9}
~& $\lambda_{\min}$ & $\lambda_{\max}$ & $\kappa$ & \emph{ratio} &
$\lambda_{\min}$ & $\lambda_{\max}$ & $\kappa$ & \multicolumn{1}{c||}{\emph{ratio}} \\\hline

 64 & 1.938E-2 & 6.982E-1 & 3.603E+1 &   -   & 3.049E-2 & 9.275E+0 & 3.042E+2 & \multicolumn{1}{c||}{  -  } \\ \hline
128 & 8.941E-3 & 5.648E-1 & 6.316E+1 & 0.570 & 1.315E-2 & 1.065E+1 & 8.101E+2 & \multicolumn{1}{c||}{0.376} \\ \hline
256 & 4.252E-3 & 4.576E-1 & 1.076E+2 & 0.587 & 5.881E-3 & 1.224E+1 & 2.081E+3 & \multicolumn{1}{c||}{0.389} \\ \hline
512 & 2.061E-3 & 3.712E-1 & 1.801E+2 & 0.598 & 2.705E-3 & 1.405E+1 & 5.196E+3 & \multicolumn{1}{c||}{0.400} \\ \hline

\end{tabular}
\end{table}

\begin{table}[htbp]
\small
\centering\caption{The smallest and largest eigenvalues, condition numbers and ratios with $\tau=h^2$.}\label{ctp-04-06}\vskip 0.1cm
\begin{tabular}{||c|c|c|c|c|c|c|c|c|c|}\hline
\multirow{2}{*}{$M$}&\multicolumn{4}{c|}{$\alpha_0=0.9$, $\alpha_1=0.4$, $\beta=0.3$, $\gamma=0.8$}
&\multicolumn{4}{c||}{$\alpha_0=0.7$, $\alpha_1=0.5$, $\beta=0.15$, $\gamma=0.95$} \\ \cline{2-9}
~& $\lambda_{\min}$ & $\lambda_{\max}$ & $\kappa$ & \emph{ratio} &
$\lambda_{\min}$ & $\lambda_{\max}$ & $\kappa$ & \multicolumn{1}{c||}{\emph{ratio}} \\\hline

 32 & 3.230E-2 & 4.853E-2 & 1.503E+0 &   -   & 4.060E-2 & 7.249E-1 & 1.785E+1 & \multicolumn{1}{c||}{  -  } \\ \hline
 64 & 1.585E-2 & 2.191E-2 & 1.382E+0 & 1.087 & 1.870E-2 & 5.103E-1 & 2.729E+1 & \multicolumn{1}{c||}{0.654} \\ \hline
128 & 7.864E-3 & 1.003E-2 & 1.275E+0 & 1.084 & 8.888E-3 & 3.596E-1 & 4.046E+1 & \multicolumn{1}{c||}{0.674} \\ \hline
256 & 3.918E-3 & 4.662E-3 & 1.190E+0 & 1.072 & 4.297E-3 & 2.537E-1 & 5.904E+1 & \multicolumn{1}{c||}{0.685} \\ \hline

\end{tabular}
\end{table}

\begin{table}[htbp]
\small
\centering\caption{The smallest and largest eigenvalues, condition numbers and ratios with $\tau=1/64$.}\label{ctp-04-07}\vskip 0.1cm
\begin{tabular}{||c|c|c|c|c|c|c|c|c|c|}\hline
\multirow{2}{*}{$M$}&\multicolumn{4}{c|}{$\alpha_0=0.9$, $\alpha_1=0.4$, $\beta=0.3$, $\gamma=0.8$}
&\multicolumn{4}{c||}{$\alpha_0=0.7$, $\alpha_1=0.5$, $\beta=0.15$, $\gamma=0.95$} \\ \cline{2-9}
~& $\lambda_{\min}$ & $\lambda_{\max}$ & $\kappa$ & \emph{ratio} &
$\lambda_{\min}$ & $\lambda_{\max}$ & $\kappa$ & \multicolumn{1}{c||}{\emph{ratio}} \\\hline

 64 & 1.938E-2 & 6.982E-1 & 3.603E+1 &   -   & 3.049E-2 & 9.275E+0 & 3.042E+2 & \multicolumn{1}{c||}{  -  } \\ \hline
128 & 9.691E-3 & 1.052E+0 & 1.085E+2 & 0.332 & 1.525E-2 & 1.730E+1 & 1.135E+3 & \multicolumn{1}{c||}{0.268} \\ \hline
256 & 4.846E-3 & 1.590E+0 & 3.281E+2 & 0.331 & 7.625E-3 & 3.229E+1 & 4.234E+3 & \multicolumn{1}{c||}{0.268} \\ \hline
512 & 2.423E-3 & 2.408E+0 & 9.939E+2 & 0.330 & 3.813E-3 & 6.025E+1 & 1.580E+4 & \multicolumn{1}{c||}{0.268} \\ \hline

\end{tabular}
\end{table}

\begin{Example}\label{cex-04-04}
  We use the same Example \ref{cex-04-02}.
\end{Example}

This example reveals the effect of $K_2$ on condition numbers of typical cases
$\tau=h$, $\tau=h^2$ and $\tau=1/64$ for $\alpha_0=0.7$, $\alpha_1=0.4$, $\beta=0.3$, $\gamma=0.85$
and $\alpha_0=0.8$, $\alpha_1=0.3$, $\beta=0.2$, $\gamma=0.75$. As expected,
the results shown in Table \ref{ctp-04-08} declare that condition numbers are independent of $K_2$ when $K_2$ is large enough.

\begin{table}[htbp]
\footnotesize
\centering\caption{The condition numbers and ratios for Example \ref{cex-04-04}.}\label{ctp-04-08}\vskip 0.1cm
\begin{tabular}{||c|c|c|c|c|c|c|c|c|c|c|c|c|}\hline
\multirow{3}{*}{$K_2$}&\multicolumn{6}{c|}{$\alpha_0=0.7$, $\alpha_1=0.4$, $\beta=0.3$, $\gamma=0.85$}
&\multicolumn{6}{c||}{$\alpha_0=0.8$, $\alpha_1=0.3$, $\beta=0.2$, $\gamma=0.75$} \\ \cline{2-13}
~&\multicolumn{2}{c|}{$h=\tau$}&\multicolumn{2}{c|}{$h=\sqrt{\tau}$}&\multicolumn{2}{c|}{$\tau=1/64$}
&\multicolumn{2}{c|}{$h=\tau$}&\multicolumn{2}{c|}{$h=\sqrt{\tau}$}&\multicolumn{2}{c||}{$\tau=1/64$} \\ \cline{2-13}
~& $\kappa$ & \emph{ratio} & $\kappa$ & \emph{ratio} & $\kappa$ & \emph{ratio}
& $\kappa$ & \emph{ratio} & $\kappa$ & \emph{ratio} & $\kappa$ & \multicolumn{1}{c||}{\emph{ratio}} \\\hline

3E2 & 5.06E3 &  -   & 3.77E2 &  -   & 5.21E3 &  -   & 1.63E3 &   -  & 1.12E2 &  -   & 1.75E3 & \multicolumn{1}{c||}{ -  } \\ \hline
3E3 & 5.31E3 & 0.95 & 4.89E2 & 0.77 & 5.33E3 & 0.98 & 1.81E3 & 0.90 & 2.07E2 & 0.54 & 1.82E3 & \multicolumn{1}{c||}{0.96} \\ \hline
3E4 & 5.34E3 & 0.99 & 5.04E2 & 0.97 & 5.34E3 & 0.99 & 1.83E3 & 0.99 & 2.26E2 & 0.91 & 1.83E3 & \multicolumn{1}{c||}{0.99} \\ \hline
3E5 & 5.34E3 & 1.00 & 5.06E2 & 1.00 & 5.34E3 & 1.00 & 1.83E3 & 1.00 & 2.28E2 & 0.99 & 1.83E3 & \multicolumn{1}{c||}{1.00} \\ \hline

\end{tabular}
\end{table}

\subsection{Performance evaluation test}

\begin{Example}\label{cex-04-05}
  Comparisons of $i$CAMG over CG and CAMG with Ruge-St{\"u}ben coarsening strategy and $\mathcal{U}_0=0$.
\end{Example}

\begin{table}[htbp]
\footnotesize
\centering\caption{Number of iterations and wall time for the case $\tau=h$.}\label{ctp-04-09}\vskip 0.1cm
\begin{tabular}{||c|c|c|c|c|c|c|c|c|c|c|c|c|c|c|c|c|}\hline
\multirow{3}{*}{$M$}&\multicolumn{6}{c|}{$\alpha_0=0.9$, $\alpha_1=0.4$, $\beta=0.3$, $\gamma=0.8$}
&\multicolumn{6}{c||}{$\alpha_0=0.7$, $\alpha_1=0.5$, $\beta=0.15$, $\gamma=0.95$} \\ \cline{2-13}
~&\multicolumn{2}{c|}{CG}&\multicolumn{2}{c|}{CAMG}&\multicolumn{2}{c|}{$i$CAMG}
 &\multicolumn{2}{c|}{CG}&\multicolumn{2}{c|}{CAMG}&\multicolumn{2}{c||}{$i$CAMG} \\ \cline{2-13}
~& \emph{Its} & $T_c$ & \emph{Its} & $T_c$ & \emph{Its} & $T_c$ & \emph{Its} & $T_c$
& \emph{Its} & $T_c$ & \emph{Its} & \multicolumn{1}{c||}{$T_c$} \\\hline

 512 & 151 & 3.23E-2 & 7 & 1.54E-2 & 7 & 5.52E-3 & 249 & 5.43E-2 & 5 & 1.31E-2 & 5 & \multicolumn{1}{c||}{4.58E-3} \\ \hline
1024 & 225 & 9.41E-2 & 7 & 6.93E-2 & 8 & 1.63E-2 & 479 & 2.08E-1 & 5 & 6.97E-2 & 5 & \multicolumn{1}{c||}{1.29E-2} \\ \hline
2048 & 300 & 2.95E-1 & 7 & 2.86E-1 & 9 & 4.79E-2 & 920 & 9.19E-1 & 5 & 2.89E-1 & 5 & \multicolumn{1}{c||}{3.75E-2} \\ \hline
4096 & 385 & 6.99E-1 & 7 & 1.27E+0 & 9 & 1.38E-1 & $\ast$ & $\ast$ & 5 & 1.21E+0 & 5 & \multicolumn{1}{c||}{1.10E-1} \\ \hline

\end{tabular}
\end{table}

\begin{table}[htbp]
\footnotesize
\centering\caption{Number of iterations and wall time for the case $\tau=h^2$.}\label{ctp-04-10}\vskip 0.1cm
\begin{tabular}{||c|c|c|c|c|c|c|c|c|c|c|c|c|}\hline
\multirow{3}{*}{$M$}&\multicolumn{6}{c|}{$\alpha_0=0.9$, $\alpha_1=0.4$, $\beta=0.3$, $\gamma=0.8$}
&\multicolumn{6}{c||}{$\alpha_0=0.7$, $\alpha_1=0.5$, $\beta=0.15$, $\gamma=0.7$} \\ \cline{2-13}
~&\multicolumn{2}{c|}{CG}&\multicolumn{2}{c|}{CAMG}&\multicolumn{2}{c|}{$i$CAMG}
 &\multicolumn{2}{c|}{CG}&\multicolumn{2}{c|}{CAMG}&\multicolumn{2}{c||}{$i$CAMG} \\ \cline{2-13}
~& \emph{Its} & $T_c$ & \emph{Its} & $T_c$ & \emph{Its} & $T_c$
& \emph{Its} & $T_c$ & \emph{Its} & $T_c$ & \emph{Its} & \multicolumn{1}{c||}{$T_c$} \\\hline

 512 & 8 & 1.62E-3 & 10 & 7.19E-3 & 8 & 1.12E-3 & 17 & 3.82E-3 & 7 & 1.39E-2 & 17 & \multicolumn{1}{c||}{3.82E-3} \\ \hline
1024 & 8 & 3.58E-3 & 8  & 3.29E-2 & 8 & 3.58E-3 & 17 & 7.01E-3 & 7 & 6.49E-2 & 17 & \multicolumn{1}{c||}{7.01E-3} \\ \hline
2048 & 8 & 8.01E-3 & 8  & 1.07E-1 & 8 & 8.81E-3 & 17 & 1.73E-2 & 7 & 2.62E-1 & 17 & \multicolumn{1}{c||}{1.73E-2} \\ \hline
4096 & 9 & 1.82E-2 & 9  & 4.70E-1 & 9 & 1.82E-2 & 17 & 3.29E-2 & 7 & 1.04E+0 & 17 & \multicolumn{1}{c||}{3.29E-2} \\ \hline

\end{tabular}
\end{table}

\begin{table}[htbp]
\footnotesize
\centering\caption{Number of iterations and wall time for the case $\tau=1/64$.}\label{ctp-04-11}\vskip 0.1cm
\begin{tabular}{||c|c|c|c|c|c|c|c|c|c|c|c|c|c|c|c|c|}\hline
\multirow{3}{*}{$M$}&\multicolumn{6}{c|}{$\alpha_0=0.9$, $\alpha_1=0.4$, $\beta=0.3$, $\gamma=0.8$}
&\multicolumn{6}{c||}{$\alpha_0=0.7$, $\alpha_1=0.5$, $\beta=0.15$, $\gamma=0.95$} \\ \cline{2-13}
~&\multicolumn{2}{c|}{CG}&\multicolumn{2}{c|}{CAMG}&\multicolumn{2}{c|}{$i$CAMG}
 &\multicolumn{2}{c|}{CG}&\multicolumn{2}{c|}{CAMG}&\multicolumn{2}{c||}{$i$CAMG} \\ \cline{2-13}
~& \emph{Its} & $T_c$ & \emph{Its} & $T_c$ & \emph{Its} & $T_c$
& \emph{Its} & $T_c$ & \emph{Its} & $T_c$ & \emph{Its} & \multicolumn{1}{c||}{$T_c$} \\\hline

 512 & 173 & 3.69E-2 & 7 & 1.65E-2 & 7 & 5.56E-3 & 250 & 5.38E-2 & 5 & 1.25E-2 & 4 & \multicolumn{1}{c||}{3.88E-3} \\ \hline
1024 & 301 & 1.26E-1 & 8 & 7.55E-2 & 7 & 1.55E-2 & 483 & 2.04E-1 & 5 & 5.65E-2 & 4 & \multicolumn{1}{c||}{1.10E-2} \\ \hline
2048 & 524 & 5.32E-1 & 8 & 3.10E-1 & 7 & 4.48E-2 & 933 & 9.37E-1 & 5 & 2.31E-1 & 4 & \multicolumn{1}{c||}{3.07E-2} \\ \hline
4096 & 908 & 1.66E+0 & 8 & 1.31E+0 & 7 & 1.25E-1 & $\ast$ & $\ast$ & 5 & 1.01E+0 & 4 & \multicolumn{1}{c||}{8.71E-2} \\ \hline

\end{tabular}
\end{table}

Table \ref{ctp-04-09}, \ref{ctp-04-10} and \ref{ctp-04-11} respectively give the results for cases $\tau=h$, $\tau=h^2$ and $\tau=1/64$,
which illustrate that CAMG and $i$CAMG both converge robustly with respect to mesh sizes and fractional orders,
while CG is only suitable for $\tau=h^2$ because of \eqref{chp-03-05}.
$i$CAMG adaptively adjusts to be CG in such circumstance.
Here our emphasis is on computational effort. The cost of $i$CAMG increases $\mathcal{O}(M\log M)$,
yet CAMG bears $\mathcal{O}(M^2)$ to achieve convergence. Hence, $i$CAMG exhibits a significant advantage over CAMG,
runs 9.2 and 11.0 times faster for $\tau=h$, 25.8 and 31.6 for $\tau=h^2$, 10.5 and 11.6 for $\tau=1/64$ at the size of $M=4096$.

\section{Conclusion}

Classical AMG is an $\mathcal{O}(M^2)$ solution process for SFDEs. We have proposed a lossless in robustness and adaptive variant
with $\mathcal{O}(M\log M)$ algorithmic complexity and $\mathcal{O}(M)$ matrix-free storage,
employed to solve the time-space FE discretization of one-dimensional MTFADEs. The approach relied on a number of
theoretically proved characterizations and condition number estimation on the resulting matrix,
an effective measure on the strength-of-connection tolerance and straightforward modifications to grid-transfer operators.
Our numerical results verify the saturation error order of the discretization,
well robustness and considerable advantages of the proposed solver over CG and classical AMG methods.

\section*{Acknowledgments}

This work is under auspices of National Natural Science Foundation of China (11571293, 11601460, 11601462)
and the General Project of Hunan Provincial Education Department of China (16C1540, 17C1527).


\section*{References}

\begin{thebibliography}{99}

\bibitem{m-001}
K. S. Miller and B. Ross,
An introduction to the fractional calculus and fractional differential equations,
Wiley (1993).

\bibitem{l-003}
F. Liu, P. Zhuang, V. Anh, I. Turner, and K. Burrage,
Stability and convergence of the difference methods for the space-time fractional advection-diffusion equation,
Appl. Math. Comput. 191 (1) (2007), pp. 12-20.

\bibitem{g-001}
G. H. Gao and Z. Z. Sun,
A compact finite difference scheme for the fractional sub-diffusion equations,
J. Comput. Phys. 230 (3) (2011), pp. 586-595.

\bibitem{c-002}
X. N. Cao, J. L. Fu and H. Huang,
Numerical method for the time fractional Fokker-Planck equation,
Adv. Appl. Math. Mech. 4 (6) (2012), pp. 848-863.

\bibitem{w-001}
H. Wang and T. S. Basu,
A fast finite difference method for two-dimensional space-fractional diffusion equations,
SIAM J. Sci. Comput. 34 (5) (2012), pp. A2444-A2458.

\bibitem{y-002}
W. Yang, D. L. Wang and L. Yang,
A stable numerical method for space fractional Landau-Lifshitz equations,
Appl. Math. Lett. 61 (2016), pp. 149-155.

\bibitem{e-001}
V. J. Ervin and J. P. Roop,
Variational formulation for the stationary fractional advection dispersion equation,
Numer. Methods Partial Differ. Equ. 22 (3) (2006), pp. 558-576.

\bibitem{z-002}
H. Zhang, F. Liu and V. Anh,
Galerkin finite element approximation of symmetric space-fractional partial differential equations,
Appl. Math. Comput. 217 (6) (2010), pp. 2534-2545.

\bibitem{b-002}
K. Burrage, N. Hale and D. Kay,
An efficient implicit FEM scheme for fractional-in-space reaction-diffusion equations,
SIAM J. Sci. Comput. 34 (4) (2012), pp. A2145-A2172.

\bibitem{b-001}
W. P. Bu, Y. F. Tang and J. Y. Yang,
Galerkin finite element method for two-dimensional Riesz space fractional diffusion equations,
J. Comput. Phys. 276 (2014), pp. 26-38.

\bibitem{m-002}
K. Mustapha, B. Abdallah and K. M. Furati,
A discontinuous Petrov-Galerkin method for time-fractional diffusion equations,
SIAM J. Numer. Anal. 52 (5) (2014), pp. 2512-2529.

\bibitem{f-001}
L. B. Feng, P. Zhuang, F. Liu, I. Turner, and Y. T. Gu,
Finite element method for space-time fractional diffusion equation,
Numer. Algor. 72 (3) (2016), pp. 749-767.

\bibitem{b-004}
W. P. Bu, A. G. Xiao and W. Zeng,
Finite difference/finite element methods for distributed-order time fractional diffusion equations,
J. Sci. Comput. 72 (1) (2017), pp. 422-441.

\bibitem{l-004}
F. Liu, P. Zhuang, I. Turner, K. Burrage, and V. Anh,
A new fractional finite volume method for solving the fractional diffusion equation,
Appl. Math. Model. 38 (15-16) (2014), pp. 3871-3878.

\bibitem{l-005}
Y. M. Lin and C. J. Xu,
Finite difference/spectral approximations for the time-fractional diffusion equation,
J. Comput. Phys. 225 (2) (2007), pp. 1533-1552.

\bibitem{z-004}
M. Zayernouri, W. R. Cao, Z. Q. Zhang, and G. E. Karniadakis,
Spectral and discontinuous spectral element methods for fractional delay equations,
SIAM J. Sci. Comput. 36 (6) (2014), pp. B904-B929.

\bibitem{y-005}
Y. Yang, Y. P. Chen and Y. Q. Huang,
Spectral-collocation method for fractional Fredholm integro-differential equations,
J. Korean Math. Soc. 51 (1) (2014), pp. 203-224.

\bibitem{y-003}
Y. Yang,
Jacobi spectral Galerkin methods for Volterra integral equations with weakly singular kernel,
Bull. Korean Math. Soc. 53 (1) (2016), pp. 247-262.

\bibitem{m-003}
R. Metzler and J. Klafter,
The random walk's guide to anomalous diffusion: a fractional dynamics approach,
Phys. Rep. 339 (2000), pp. 1-77.

\bibitem{s-002}
S. Shen, F. Liu, V. Anh, and I. Turner,
The fundamental solution and numerical solution of the Riesz fractional advection-dispersion equation,
IMA J. Appl. Math. 73 (6) (2008), pp. 850-872.

\bibitem{z-005}
Y. Y. Zheng, C. P. Li and Z. G. Zhao,
A note on the finite element method for the space-fractional advection diffusion equation,
Comput. Math. Appl. 59 (5) (2010), pp. 1718-1726.

\bibitem{m-004}
M. M. Meerschaert and C. Tadjeran,
Finite difference approximations for fractional advection-dispersion flow equations,
J. Comput. Appl. Math. 172 (1) (2004), pp. 65-77.

\bibitem{b-005}
A. H. Bhrawy and D. Baleanu,
A spectral Legendre-Gauss-Lobatto collocation method for a space-fractional advection diffusion equations with variable coefficients,
Rep. Math. Phys. 72 (2) (2013), pp. 219-233.

\bibitem{j-001}
H. Jiang, F. Liu, I. Turner, and K. Burrage,
Analytical solutions for the multi-term time-space Caputo-Riesz fractional advection-diffusion equations on a finite domain,
J. Math. Anal. Appl. 389 (2) (2012), pp. 1117-1127.

\bibitem{y-006}
H. Ye, F. Liu, V. Anh, and I. Turner,
Maximum principle and numerical method for the multi-term time-space Riesz-Caputo fractional differential equations,
Appl. Math. Comput. 227 (2014), pp. 531-540.

\bibitem{b-003}
W. P. Bu, X. T. Liu, Y. F. Tang, and J. Y. Yang,
Finite element multigrid method for multi-term time fractional advection diffusion equations,
Int. J. Model. Simul. Sci. Comput. 6 (2015), 1540001.

\bibitem{g-002}
C. Y. Gong, W. M. Bao, G. J. Tang, Y. W. Jiang, and J. Liu,
Computational challenge of fractional differential equations and the potential solutions: a survey,
Math. Probl. Eng. 2015 (2015), 258265.

\bibitem{t-001}
U. Trottenberg, C. W. Oosterlee and A. Sch{\"u}ller,
Multigrid,
Elsevier (2001).

\bibitem{f-002}
C. S. Feng, S. Shu, J. Xu, and C. S. Zhang,
Numerical study of geometric multigrid methods on CPU-GPU heterogeneous computers,
Adv. Appl. Math. Mech. 6 (1) (2014), pp. 1-23.

\bibitem{f-003}
R. D. Falgout, S. Friedhoff, Tz. V. Kolev, S. P. Maclachlan, and J. B. Schroder,
Parallel time integration with multigrid,
SIAM J. Sci. Comput. 36 (6) (2014), pp. C635-C661.

\bibitem{x-001}
J. Xu and L. Zikatanov,
Algebraic multigrid methods,
Acta Numer. 26 (2017), pp. 591-721.

\bibitem{p-001}
H. K. Pang and H. W. Sun,
Multigrid method for fractional diffusion equations,
J. Comput. Phys. 231 (2012), pp. 693-703.

\bibitem{z-003}
Z. Q. Zhou and H. Y. Wu,
Finite element multigrid method for the boundary value problem of fractional advection dispersion equation,
J. Appl. Math. 2013 (2013), 385463.

\bibitem{j-002}
Y. J. Jiang and X. J. Xu,
Multigrid methods for space fractional partial differential equations,
J. Comput. Phys. 302 (2015), pp. 374-392.

\bibitem{c-003}
L. Chen, R. H. Nochetto, E. Ot{\'a}rola, and A. J. Salgado,
Multilevel methods for nonuniformly elliptic operators and fractional diffusion,
Math. Comput. 85 (302) (2016), pp. 2583-2607.

\bibitem{z-006}
X. Zhao, X. Z. Hu, W. Cai, and G. E. Karniadakis,
Adaptive finite element method for fractional differential equations using hierarchical matrices,
Comput. Methods Appl. Mech. Engrg. 325 (2017), pp. 56-76.

\bibitem{c-001}
M. H. Chen and W. H. Deng,
Convergence proof for the multigrid method of the nonlocal model,
arXiv: 1605.05481.

\bibitem{y-001}
X. Q. Yue, W. P. Bu, S. Shu, M. H. Liu, and S. Wang,
Fully finite element adaptive algebraic multigrid method for time-space Caputo-Riesz fractional diffusion equations,
arXiv: 1707.08345.

\bibitem{s-003}
Y. Saad,
Iterative methods for sparse linear systems,
SIAM (2003).

\bibitem{x-002}
J. Xu and L. Zikatanov,
Algebraic multigrid methods,
Acta Numer. 26 (2017), pp. 591-721.

\bibitem{h-001}
V. E. Henson and U. M. Yang,
BoomerAMG: A parallel algebraic multigrid solver and preconditioner,
Appl. Numer. Math. 41 (2002), pp. 155-177.

\bibitem{m-005}
M. Naumov et al.,
AmgX: A library for GPU accelerated algebraic multigrid and preconditioned iterative methods,
SIAM J. Sci. Comput. 37 (5) (2015), pp. S602-S626.

\end{thebibliography}
\end{document}